\documentclass[reqno,a4paper,10pt]{article}

\usepackage[hmarginratio=1:1, top=2.5cm, bottom=2.5cm, left=3cm, right=3cm]{geometry}

\usepackage{amsmath}
\usepackage{amsthm} 
\usepackage{amssymb} 
\usepackage{amsfonts}
\usepackage{graphicx} 
\usepackage[english]{babel}
\usepackage{enumerate}
\usepackage{listings} 
\usepackage{bm}
\usepackage{float}
\usepackage{hyperref}
\usepackage{bbm}
\usepackage{tikz}
\usepackage{nicematrix}
\usepackage{authblk}

\tikzstyle{vertex}=[circle, draw, fill=black, inner sep=0pt, minimum size=4pt]
\tikzstyle{ivertex}=[circle, draw, fill=black, inner sep=0pt, minimum size=0pt]
\tikzstyle{blankvertex}=[circle, draw=white, fill=white, inner sep=0pt, minimum size=4pt]
\tikzstyle{edge}=[line width=1.5pt]
\tikzstyle{iedge}=[line width=0pt]
\tikzstyle{labelsty}=[font=\scriptsize]
\tikzstyle{dedge}=[edge,-latex]
\tikzstyle{pedge}=[dashed width=1.5pt]

\usepackage{mathrsfs}
\usepackage{blkarray}
\usepackage{caption}
\usepackage{subcaption}
\usepackage{booktabs}

\newtheorem{theorem}{Theorem}[section]
\newtheorem{lemma}[theorem]{Lemma}

\newtheorem{definition}[theorem]{Definition}

\newtheorem{proposition}[theorem]{Proposition}
\theoremstyle{definition}

\theoremstyle{remark}

\numberwithin{equation}{section}
\numberwithin{figure}{section}
\numberwithin{table}{section}
\AtBeginDocument{\numberwithin{lstlisting}{section}}

\newcommand{\R}{\mathbb{R}}

\newcommand{\Z}{\mathbb{Z}}
\newcommand{\N}{\mathbb{N}}

\newcommand{\G}{\tilde{G}}
\newcommand{\p}{\tilde{p}}

\title{Non-Euclidean Crystallographic Rigidity}

\author[1,2]{Jack Esson}
\affil[1]{School of Mathematical Sciences, Lancaster University, Lancaster LA1 4YF, United Kingdom.}
\affil[2]{Email address: j.esson@lancaster.ac.uk (corresponding author)}

\author[1,3]{Eleftherios Kastis}
\affil[3]{Email address: l.kastis@lancaster.ac.uk}

\author[1,4]{Bernd Schulze}
\affil[4]{Email address: b.schulze@lancaster.ac.uk}

\begin{document}

	\raggedbottom

    \maketitle

    \begin{abstract}
        This paper establishes combinatorial characterisations of forced-symmetric and forced-periodic rigidity (under a fixed lattice) of bar-joint frameworks in non-Euclidean normed planes. In $\ell_q$-planes for $q\in(1,\infty)\backslash\{2\}$, we prove characterisations for forced-periodic rigidity and forced-reflectionally-symmetric rigidity. We also characterise forced-symmetric rigidity in this space with respect to the orientation-reversing wallpaper group $\Z^2\rtimes\mathcal{C}_s$, otherwise known as $pm$ in crystallography. In the $\ell_1$ and $\ell_\infty$-planes, we provide characterisations for forced-periodic rigidity and forced-$\Z^2\rtimes\mathcal{C}_s$-symmetric rigidity. All of these characterisations are proved by inductive constructions involving Henneberg-type graph operations.
    \end{abstract}


    \section{Introduction and Preliminaries}

    \subsection{Introduction}\label{SectionIntroduction}
    A $2$-dimensional \emph{(bar-joint) framework} is a pair $(G,p)$ consisting of a simple graph $G=(V(G),E(G))$ and a map $p:V(G)\to \mathbb{R}^2$ assigning positions to the vertices of $G$, with $p(v_i)\neq p(v_j)$ for $\{v_i,v_j\}\in E(G)$. The map $p$ is also called a \emph{configuration} of $G$ in $\mathbb{R}^2$. Bar-joint frameworks can be used to model real-world structures made of stiff bars connected by freely rotational joints. These appear in a variety of applications, including in engineering, materials science, structural biology and crystallography. A primary interest in this area is to characterise when a bar-joint framework is rigid, meaning that it cannot be continuously deformed while maintaining the edge lengths. For a summary of the basics of rigidity theory in Euclidean spaces, see 
\cite{HandbookDCG,connelly_guest_2022,WhiteleyMatroids}.

    Over the last two decades, there has been significant interest in the forced-symmetric rigidity of symmetric frameworks, including the forced-periodic rigidity of periodic frameworks. For forced-symmetric rigidity, one only considers motions of the framework that preserve the original symmetry. A separate area of study is incidental-symmetric rigidity, where motions are not required to preserve the symmetry, but this is beyond the scope of this paper. When the configurations are assumed to be as generic as possible within the symmetry constraints, forced-symmetric rigidity in Euclidean spaces can be characterised in terms of  group-labelled quotient graphs. For a summary of combinatorial results regarding the  rigidity of finite symmetric frameworks in Euclidean spaces, we refer the reader to \cite{HandbookDCGsym,constraint19}.
    
    For infinite periodic frameworks,  the rigidity analysis has largely focused on the forced-periodic setting. Here the theory splits into two parts, depending on whether the lattice representation is fixed or is allowed to be flexible as the framework moves. In this article, we will focus only on the case where the lattice is fixed. For periodic frameworks on a fixed lattice in the Euclidean plane, a combinatorial characterisation for minimal forced-periodic rigidity was given by E. Ross in \cite{Inductive}. In a previous paper, we characterised conditions for forced-symmetric rigidity in the Euclidean plane with respect to the orientation-reversing wallpaper group $\Z^2\rtimes\mathcal{C}_s$ \cite{EssonCrystallographic}. Characterisations for forced-symmetric rigidity in the Euclidean plane with respect to orientation-preserving wallpaper groups involving rotations were covered by D. Bernstein in \cite{BernsteinRotations}.

    In this paper, we are interested in finding similar results regarding forced-symmetric and forced-periodic rigidity for non-Euclidean $\ell_q$-planes, including the $\ell_1$-plane and the $\ell_\infty$-plane. The rigidity of finite non-symmetric frameworks in these spaces has received a lot of attention lately, see e.g. \cite{NonEuclidean,KitsonDenseOpen,Dewar,dhn24}. Note that we are only considering $2$-dimensional spaces here. Currently, no combinatorial characterisation is known for finite non-symmetric rigidity for bar-joint frameworks in any higher-dimensional $\ell_q$-spaces (for $q\in[1,\infty]$), so there are no characterisations for symmetric or periodic rigidity either. (Note that there exist results for rigidity in higher-dimensional spaces with mixed norms, see \cite{DewarCylindrical}. There also exist results for the special class of \emph{body-bar} frameworks in all dimensions, both for finite symmetric ones and for periodic ones under a fixed lattice representation, see \cite{TanigawaBodyBar}.)
    
    Like in Euclidean spaces, it is easier to linearise the problem of rigidity by differentiating the length constraints and considering infinitesimal rigidity. The theory of infinitesimal rigidity in general normed spaces is described in \cite{NonEuclidean} and we briefly summarise it here. In general, an \emph{infinitesimal motion} of a framework $(G,p)$ in a normed plane $(\R^2,\|\cdot\|)$ is a map $u:V(G)\to\R^2$ such that, for all $\{v_i,v_j\}\in E(G)$, as $t\to0$,
    \begin{align*}
        \|(p(v_i)+tu(v_i))-(p(v_j)+tu(v_j))\|-\|p(v_i)-p(v_j)\|\to o(t).
    \end{align*}
    An infinitesimal motion is \emph{trivial} if it corresponds to a rigid  motion of $(\R^2,\|\cdot\|)$ (i.e. a family 
    $\alpha=\{\alpha_x\}_{x\in\mathbb{R}^2}$ 
    of continuous paths, differentiable at time 0 and satisfying $\alpha_x(0)=x$ for all $x\in \mathbb{R}^2$, so that the distance between any pair of points remains the same along the paths) 
    and \emph{non-trivial} otherwise \cite[Definition 2.2]{NonEuclidean}. A framework is \emph{infinitesimally rigid} if all of its infinitesimal motions are trivial and \emph{infinitesimally flexible} otherwise.

    Our settings for this paper will be planes with $\ell_q$-norms for $q\in(1,\infty)\backslash\{2\}$ and with the $\ell_\infty$ (or equivalently, $\ell_1$) norm. With all of these norms, the space of trivial infinitesimal motions is $2$-dimensional, with a basis consisting of two translations \cite{NonEuclidean}. This is in contrast to the Euclidean plane, where there is a $3$-dimensional space of trivial infinitesimal motions, as rotations also yield trivial infinitesimal motions.
    
    Note that  there is no notion of ``genericity" of a configuration in the $\ell_1$-plane or $\ell_\infty$-plane; there may be open sets of configurations where the framework is rigid and other open sets where the framework is flexible. In other words, the set of regular configurations (where the rigidity matrix has maximum rank) of a particular graph in these spaces need not be dense in the space of all configurations. However, the set of regular configurations will always be open. In contrast, in $\ell_q$-planes for $q\in(1,\infty)$, the set of regular configurations is always dense in the space of all configurations, so the usual ideas of genericity from Euclidean rigidity can be applied to these spaces \cite{KitsonDenseOpen}.

    It is also important to note that the only relevant types of symmetry in a given normed space are those that are compatible with the symmetries of the unit sphere. In the Euclidean plane, a characterisation of forced-symmetric rigidity for any finite-order rotational symmetry was found by J. Malestein and L. Theran in \cite[Theorem 1]{MalesteinSymmetry}. However, only $2$-fold and $4$-fold rotational symmetry can be considered in non-Euclidean $\ell_q$-planes (including $\ell_1$ and $\ell_\infty$). Likewise, symmetry with respect to dihedral groups of order $2k$, where $k$ is odd, is not relevant in these planes, although a characterisation of forced-symmetric rigidity for these groups in the Euclidean plane was found by T. Jord\'{a}n, V. Kaszanitzky and S. Tanigawa in \cite[Theorem 8.2]{EGRES}. No characterisations for forced-symmetric rigidity in the Euclidean plane are known for dihedral groups of order $2k$, where $k$ is even.
    
    In \cite[Theorem 3.6]{NonEuclidean}, D. Kitson and S. Power proved that a regular framework in a non-Euclidean $\ell_q$-plane is minimally rigid if and only if its graph is $(2,2)$-tight. They also showed that the same condition is necessary and sufficient for minimal rigidity in polytopic planes with a quadrilateral unit circle, including the $\ell_1$-plane and the $\ell_\infty$-plane \cite[Theorem 4.10]{NonEuclidean}. This work opened questions regarding forced-symmetric rigidity in these non-Euclidean planes. Forced-reflectionally-symmetric rigidity in polytopic planes with a quadrilateral unit circle was characterised by D. Kitson and B. Schulze in \cite[Theorem 22]{GridLike}. A characterisation for half-turn forced-rotationally-symmetric rigidity in polytopic planes with a quadrilateral unit circle was proved by D. Kitson, A. Nixon and B. Schulze in \cite[Theorem 4.3]{KitsonSymmetricNormed}.

In this paper,  we prove characterisations for forced-periodic rigidity (Theorem \ref{lqPeriodic}), forced-reflectionally-symmetric rigidity (Theorem \ref{Reflectivelq}) and forced-$\Z^2\rtimes\mathcal{C}_s$-symmetric rigidity (Theorem \ref{ZReflectivelq}) in $\ell_q$-planes with $q\in(1,\infty)\backslash\{2\}$. For the $\ell_1$ and $\ell_\infty$-planes, we prove characterisations for forced-periodic rigidity (Theorem \ref{PolyPeriodic}) and forced-$\Z^2\rtimes\mathcal{C}_s$-symmetric rigidity (Theorem \ref{ZReflectivePoly}). Table \ref{TableSymmetry} summarises known and unknown results for forced-symmetric rigidity in these spaces with respect to various symmetry groups.

    \begin{table}[H]
        \begin{center}
        \begin{tabular}{|l||l|l|l|}
        \hline
        Symmetry&$\ell_2$-plane&$\ell_q$-plane for $q\in(1,\infty)\backslash\{2\}$&$\ell_1$-plane and $\ell_\infty$-plane\\
        \hline
        No symmetry&\cite{Pollaczek-Geiringer}\cite{Laman}&\cite[Theorem 3.6]{NonEuclidean}&\cite[Theorem 4.10]{NonEuclidean}\\
        $\mathcal{C}_s$&\cite[Theorem 2]{MalesteinSymmetry}&Theorem \ref{Reflectivelq}&\cite[Theorem 22]{GridLike}\\
        $\mathcal{C}_2$&\cite[Theorem 1]{MalesteinSymmetry}&Unknown&\cite[Theorem 4.3]{KitsonSymmetricNormed}\\
        $\mathcal{C}_4$&\cite[Theorem 1]{MalesteinSymmetry}&Unknown&Unknown\\
        $\Z^2$&\cite[Theorem 5.1]{Inductive}&Theorem \ref{lqPeriodic}&Theorem \ref{PolyPeriodic}\\
        $\Z^2\rtimes\mathcal{C}_s$&\cite[Theorem 3.2]{EssonCrystallographic}&Theorem \ref{ZReflectivelq}&Theorem \ref{ZReflectivePoly}\\
        $\Z^2\rtimes\mathcal{C}_2$&\cite[Theorem 5.10]{BernsteinRotations}&Unknown&Unknown\\
        $\Z^2\rtimes\mathcal{C}_4$&\cite[Theorem 5.10]{BernsteinRotations}&Unknown&Unknown\\
        \hline
    \end{tabular}
    \end{center}
    \caption{A table showing known and unknown characterisations of forced-symmetric rigidity for different symmetry groups in  normed planes.}
    \label{TableSymmetry}
    \end{table}

    Our characterisations are  variations of standard inductive arguments, including Henneberg-type graph extension moves, which are used to construct all gain graphs that satisfy the relevant conditions. To this end, we introduce new orbit rigidity matrices for symmetric frameworks in the relevant spaces. We then need to prove that the relevant Henneberg-type extension moves preserve minimal rigidity. Some of these geometric proofs require new techniques, building on work in \cite{NonEuclidean,EGRES,Surfaces, Inductive}. For the $\ell_1$ and $\ell_\infty$-planes we use symmetric variants of the framework colouring technique seen in \cite{KitsonSymmetricNormed, NonEuclidean,GridLike}. Some of these proofs require new inductive constructions and reduction arguments, specifically the proofs for $\Z^2\rtimes\mathcal{C}_s$ in the $\ell_q$-planes with $q\in(1,\infty)\backslash\{2\}$ (Section \ref{SectionLqWallpaper}) and for $\Z^2$ and $\Z^2\rtimes\mathcal{C}_s$ in the $\ell_1$ and $\ell_\infty$-planes (Sections \ref{SectionPolytopicPeriodic} and \ref{SectionPolytopicWallpaper} respectively). 

   The paper is organised as follows. In Section \ref{SectionSymmetric}, we introduce the basic theory of symmetric frameworks. Section \ref{SectionExtensions} introduces the basic extension moves, which form a larger gain graph from a smaller one. These are used for the inductive constructions throughout the paper. In Section \ref{SectionLqBasics}, we introduce the basic theory needed to understand rigidity and symmetric rigidity in $\ell_q$-planes for $q\in(1,\infty)\backslash\{2\}$. In Section \ref{SectionLqPeriodic}, this is used to characterise conditions for periodic rigidity in these spaces. Likewise, Sections \ref{SectionLqReflection} and \ref{SectionLqWallpaper} characterise conditions for $\mathcal{C}_s$-symmetric rigidity and $\Z^2\rtimes\mathcal{C}_s$-symmetric rigidity respectively in these $\ell_q$-planes. Section \ref{SectionPolytopicBasics} introduces the basics of rigidity in the $\ell_1$-plane and the $\ell_\infty$-plane. By isometric isomorphism, this theory also applies to other polytopic planes with two basis vectors. We build on this in Section \ref{SectionPolytopicPeriodic} to characterise periodic rigidity in the $\ell_1$-plane and $\ell_\infty$-plane and then in Section \ref{SectionPolytopicWallpaper} to characterise $\Z^2\rtimes\mathcal{C}_s$-symmetric rigidity in $\ell_1$-plane and $\ell_\infty$-plane. Finally, Section \ref{SectionFurther} discusses avenues for future research.

\subsection{Background on symmetric frameworks}\label{SectionSymmetric}

    We fix an arbitrary norm in the plane $\R^2$. Let $\Gamma$ be a group of isometries in the plane with respect to this norm. We refer to such a group as a \emph{symmetry group}. A simple graph $\G$ is \emph{$\Gamma$-symmetric} if there is a group action $\Gamma\to\mathrm{Aut}(\G)$. If $\Gamma$ contains translations, then we say that a $\Gamma$-symmetric graph is \emph{periodic} and the term ``periodic" is often used in place of ``symmetric" for these graphs. In this paper, we are only considering actions that are \emph{free} on the vertex set, meaning that no non-identity elements of $\Gamma$ fix any vertices. We are also assuming that symmetric graphs have a finite number of vertex orbits and a finite number of edge orbits, so a $\Gamma$-symmetric graph is infinite if and only if $\Gamma$ is infinite. A useful tool for studying symmetric graphs is the (group-labelled) \emph{gain graph}, which we now define. For more details, see \cite{EGRES,Torus,Quotient} for example.

    Let $G$ be a multigraph. Choose an arbitrary orientation for $G$ and let $\vec{E}(G)$ be the resulting oriented edge set. A \emph{$\Gamma$-gain assignment} for $G$ is a map $m:\vec{E}(G)\to\Gamma$ such that loops receive non-identity gains, parallel edges with the same orientation receive different gains and parallel edges with opposite orientations receive gains that are not inverses of each other. A \emph{$\Gamma$-gain graph} is a pair $(G,m)$ consisting of a multigraph $G$ and a $\Gamma$-gain assignment $m:\vec{E}(G)\to\Gamma$. For a given edge $e\in \vec{E}(G)$, the element $m(e)$ is known as the \emph{gain} of $e$. An edge $e$ in a gain graph from $v_i$ to $v_j$ with gain $m(e)$ is denoted by $e=(v_i,v_j;m(e))\in E(G)$.

    A gain graph can be used to represent a larger symmetric simple graph  by using only one vertex for each vertex orbit and one edge for each edge orbit. From a $\Gamma$-symmetric graph $\G$, we can construct its \emph{quotient $\Gamma$-gain graph} $(G,m)$, as follows. Choose any set of vertex orbit representatives to form $V(G)=\{v_1,...,v_n\}$. Each edge in $\G$ takes the form $\{\gamma_iv_i,\gamma_jv_j\}$ for some $v_i,v_j\in V(G)$ and $\gamma_i,\gamma_j\in\Gamma$. For each such edge, represent the corresponding edge orbit by adding $(v_i,v_j;{\gamma_i}^{-1}\gamma_j)$ to the gain graph. In this construction, the choice of orientation for each gained edge is unimportant, as two gain graphs can be regarded as equivalent if one can be obtained from the other by reversing edge orientations and replacing associated gains with their inverses.

    Conversely, given a $\Gamma$-gain graph $(G,m)$, we can construct its \emph{derived $\Gamma$-symmetric graph} $\G$ as follows. Let $V(\G)=\{(v,\gamma):v\in V(G),\gamma\in\Gamma\}$. For brevity, $(v,\gamma)$ is often abbreviated to $\gamma v$. Then $\{\gamma_iv_i,\gamma_jv_j\}$ is in $E(\G)$ if and only if there exists $(v_i,v_j;{\gamma_i}^{-1}\gamma_j)\in E(G)$.
    
    A \emph{$\Gamma$-symmetric configuration} of a $\Gamma$-symmetric graph $\G$ is a configuration $\p:V(\G)\to\R^2$ such that, for all $v\in V(\G)$ and $\gamma\in\Gamma$, we have $\gamma\p(v)=\p(\gamma v)$. A \emph{$\Gamma$-symmetric framework} is a pair $(\G,\p)$ consisting of a $\Gamma$-symmetric graph $\G$ and a $\Gamma$-symmetric configuration $\p$ of $\G$. Note that a configuration $p:V(G)\to\R^d$ can also be defined on a $\Gamma$-gain graph $(G,m)$. This configuration $p$ derives a $\Gamma$-symmetric configuration $\p:V(\G)\to\R^d$ of the derived $\Gamma$-symmetric graph $\G$. Conversely,  a $\Gamma$-symmetric configuration $\p:V(\G)\to\R^d$ of $\G$  naturally  restricts to a  configuration $p:V(G)\to\R^d$ on the quotient $\Gamma$-gain graph $(G,m)$.

    Given an isometry $\gamma\in\Gamma$, let $\gamma_l$ denote the linear part of $\gamma$. For a $\Gamma$-symmetric framework $(\G,\p)$ an infinitesimal motion $\tilde{u}:V(\G)\to\R^2$ of $(\G,\p)$ is said to be \emph{$\Gamma$-symmetric} if, for all $v\in V(\G)$ and $\gamma\in\Gamma$, $\tilde{u}(\gamma v)=\gamma_l\tilde{u}(v)$. A $\Gamma$-symmetric framework is \emph{$\Gamma$-symmetrically infinitesimally rigid} if it has no non-trivial $\Gamma$-symmetric infinitesimal motions. Otherwise, it is \emph{$\Gamma$-symmetrically infinitesimally flexible}. It is \emph{minimally $\Gamma$-symmetrically infinitesimally rigid} if either it is $\Gamma$-symmetrically infinitesimally rigid and has no bars, or    
    it is $\Gamma$-symmetrically infinitesimally rigid and the removal of any orbit of bars results in a $\Gamma$-symmetrically infinitesimally flexible framework. 

    If $\Gamma$ contains translations, then a $\Gamma$-symmetric framework is infinite and periodic. In this paper, we are assuming that the periodicity lattice remains fixed under all infinitesimal motions, as the isometry group $\Gamma$ is fixed.

    The \emph{orbit rigidity matrix} of a symmetric framework is a matrix that encodes the constraints that edge orbits impose on symmetric infinitesimal motions. For Euclidean spaces, the orbit rigidity matrix is defined for finite symmetry groups in \cite{Orbit} and for translation groups in \cite{Inductive}. We will formally define the orbit rigidity matrix for non-Euclidean $\ell_q$-planes with $q\in(1,\infty)\backslash\{2\}$ in Section \ref{SectionLqBasics}. The formal definition of the orbit rigidity matrix in the $\ell_1$ and $\ell_\infty$-planes will be given in Section \ref{SectionPolytopicBasics}. A more general formal definition is given in \cite{KastisSymbol}. 
    
    A $\Gamma$-symmetric configuration is said to be \emph{$\Gamma$-regular} if it achieves the maximum possible rank of the orbit rigidity matrix. Equivalently, a symmetric configuration is \emph{$\Gamma$-regular} if it achieves the minimum possible dimension of the infinitesimal motion space among all configurations of this gain graph. In $\ell_q$-normed planes for $q\in(1,\infty)$, the set of regular configurations is dense in the space of all configurations. However, this is not the case in the $\ell_1$ or $\ell_\infty$-planes. As a result, there is no notion of a gain graph being generically rigid in these spaces.

    A walk in a gain graph $(G,m)$ can be denoted by ${e_1}^{\alpha_1}{e_2}^{\alpha_2}...{e_k}^{\alpha_k}$, where $e_1,e_2,...,e_k\in E(G)$ and the values $\alpha_1,\alpha_2,...,\alpha_k$ are $+1$ for edges that are traversed forwards, and $-1$ for edges that are traversed backwards. The \emph{net gain} of such a walk is the element $\prod_{i=1}^{k}m(e_i)^{\alpha_i}$. Given some $\Gamma$-gain graph $(G,m)$ with $v\in V(G)$, the \emph{gain space} $\langle(G,m)\rangle_v$ is the subgroup of $\Gamma$ that is generated by the net gains of all closed walks in $(G,m)$ that start and end at $v$.

    A gain graph $(G,m)$ is said to be \emph{balanced} if, for all starting vertices $v\in V(G)$, the gain space $\langle(G,m)\rangle_v$ is trivial. Otherwise, the gain graph is said to be \emph{unbalanced}. The gain graph is said to be \emph{purely periodic} if, for all $v\in V(G)$, the gain space $\langle(G,m)\rangle_v$ consists only of translations.

    Let $(G,m)$ be a $\Gamma$-gain graph and choose a vertex $v\in V(G)$. Let $\gamma\in\Gamma$. A \emph{switching operation} at $v$ by $\gamma$ defines a new gain assignment $m':\vec{E}(G)\to\Gamma$ on $G$ as follows:
    \begin{align*}
			m'(e)=
			\begin{cases}
				\gamma m(e)\gamma^{-1}&\text{ if $e$ is a loop incident to $v$;}\\
				\gamma m(e)&\text{ if $e$ is a non-loop edge directed from $v$;}\\
				m(e)\gamma^{-1}&\text{ if $e$ is a non-loop edge directed to $v$;}\\
				m(e)&\text{ otherwise.}
			\end{cases}
		\end{align*}
        Intuitively, switching operations change the choice of vertex orbit representatives used in the construction of the quotient gain graph. Since the derived graph does not change, switching operations preserve rigidity and flexibility properties. They also preserve the balanced and purely periodic properties. Two gain graphs are said to be \emph{equivalent} if one can be reached from the other by a sequence of switching operations or by reversing edge orientations and inverting the corresponding gains.

        It is possible to use switching operations to assign identity gain to the edges of any spanning tree in a connected $\Gamma$-gain graph. From this, we have the following result, which is a straightforward extension of \cite[Lemma 2.3]{EGRES}.
    \begin{lemma}\label{lemma:swicth}
        Let $(G,m)$ be a $\Gamma$-gain graph such that, for every vertex $v\in V(G)$, the gain space $\langle(G,m)\rangle_v$ is contained in some subgroup $\Gamma'\leq \Gamma$. Then there is an equivalent gain graph $(G,m')$ in which every gain is an element of $\Gamma'$.
    \end{lemma}

        Let $k\in\N$ and $l,n\in\N_0$. A multigraph $G$ is said to be \emph{$(k,l)$-sparse} if all subgraphs $G'\subseteq G$ with $|E(G')|\geq1$ satisfy $|E(G')|\leq k|V(G')|-l$. The multigraph $G$ is said to be \emph{$(k,l)$-tight} if it is $(k,l)$-sparse and it satisfies $|E(G)|=k|V(G)|-l$. A gain graph $(G,m)$ is said to be \emph{$(k,l,n)$-gain tight} if it is $(k,n)$-tight and every balanced subgraph is $(k,l)$-sparse.

    \subsection{Extensions}\label{SectionExtensions}
  In this paper, we will be making use of several extension moves, which are operations that form a larger gain graph from a smaller one. This section describes the types of extension moves that will be used. These are variations of some well-known moves in symmetric rigidity theory, as seen in papers such as \cite{EGRES}, \cite{KitsonSymmetricNormed}, \cite{GridLike} and \cite{Inductive} (see also \cite{Quotient}).
  
	\begin{definition}
		Let $\Gamma$ be a symmetry group of a given $2$-dimensional normed vector space and let $(G,m)$ be a $\Gamma$-gain graph. A \emph{gained 0-extension} forms a new gain graph $(G',m')$ by letting $m'|_{E(G)}=m$ and adding a new vertex $v_0$ with incident edges $e_1=(v_0,v_1;m'(e_1))$ and $e_2=(v_0,v_2;m'(e_2))$, for some (not necessarily distinct) $v_1,v_2\in V(G)$. This move is illustrated in Figure \ref{Fig0Extension}. 
        A \emph{gained $0$-reduction} is the inverse move of a gained $0$-extension.
		
		\begin{figure}[H]
			\centering
            \begin{subfigure}{0.4\textwidth}
				\centering
				\begin{tikzpicture}[scale=0.6]
					\node[vertex,label=below:$v_1$] (v1) at (0,-2) {};
					\node[vertex,label=below:$v_2$] (v2) at (2,-2) {};
				\end{tikzpicture}
				{\Huge{$^{\mapsto}$}}
				\begin{tikzpicture}[scale=0.6]
					
					\node[vertex,label=above:$v_0$] (v0b) at (2,2) {};
					\node[vertex,label=below:$v_1$] (v1b) at (0,0) {};
					\node[vertex,label=below:$v_2$] (v2b) at (4,0) {};
					
					\draw[dedge] (v0b)edge node[left,labelsty]{$e_1$}(v1b);
					\draw[dedge] (v0b)edge node[right,labelsty]{$e_2$}(v2b);
				\end{tikzpicture}
				\caption{Two neighbours.}
			\end{subfigure}
			\begin{subfigure}{0.4\textwidth}
				\centering
				\begin{tikzpicture}[scale=0.6]
					\node[vertex,label=below:$v_1$] (v1) at (0,-2) {};
				\end{tikzpicture}
				{\Huge{$^{\mapsto}$}}
				\begin{tikzpicture}[scale=0.6]
					
					\node[vertex,label=above:$v_0$] (v0b) at (1,2) {};
					\node[vertex,label=below:$v_1$] (v1b) at (1,0) {};
					\draw[dedge] (v0b)edge [bend right=40] node[left,labelsty]{$e_1$}(v1b);
					\draw[dedge] (v0b)edge [bend left=40] node[right,labelsty]{$e_2$}(v1b);
				\end{tikzpicture}
				\caption{One neighbour.}
			\end{subfigure}
			\caption{The two variations of the gained $0$-extension.}
			\label{Fig0Extension}
		\end{figure}

A \emph{gained 1-extension} forms a new gain graph $(G',m')$ by letting $m'|_{E(G)}=m$, removing an edge $e=(v_1,v_2;m(e))\in E(G)$ and adding a new vertex $v_0$ with incident edges $e_1=(v_0,v_1;m'(e_1))$, $e_2=(v_0,v_2;m'(e_2))$ and $e_3=(v_0,v_3;m'(e_3))$, for some $v_3\in V(G)$, with the additional requirement that $(m'(e_1))^{-1}m'(e_2)=m(e)$. Note that the vertices $v_1,v_2,v_3\in V(G)$ need not be distinct. This move is illustrated in Figure \ref{Fig1Extension}. A \emph{gained $1$-reduction} is the inverse move of a gained $1$-extension.
		
		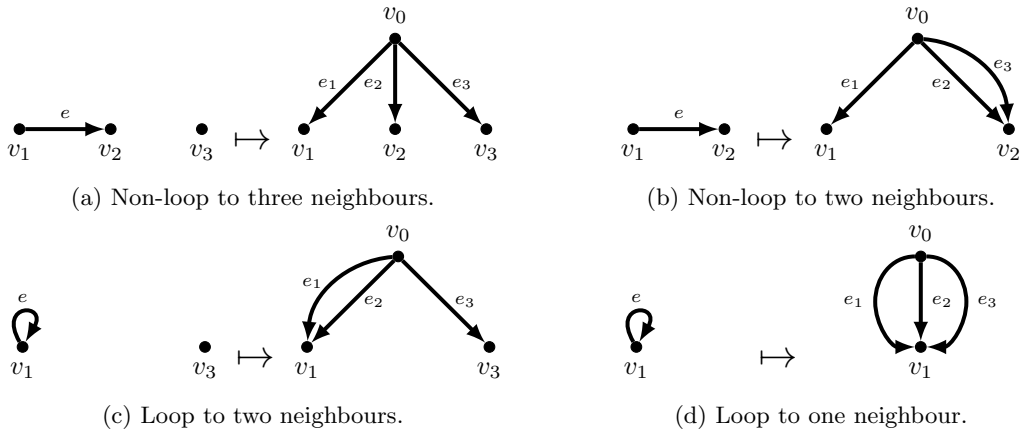
\begin{figure}[H]
			\centering

			\begin{subfigure}{0.49\textwidth}
				\centering
				\begin{tikzpicture}[scale=0.6]
					\node[vertex,label=below:$v_1$] (v1) at (-2,-1) {};
					\node[vertex,label=below:$v_2$] (v2) at (0,-1) {};
					\node[vertex,label=below:$v_3$] (v3) at (2,-1) {};
					
					\draw[dedge] (v1)edge node[above,labelsty]{$e$}(v2);
				\end{tikzpicture}
				{\Huge{$^{\mapsto}$}}
				\begin{tikzpicture}[scale=0.6]
					
					\node[vertex,label=above:$v_0$] (v0) at (0,1) {};
					\node[vertex,label=below:$v_1$] (v1) at (-2,-1) {};
					\node[vertex,label=below:$v_2$] (v2) at (0,-1) {};
					\node[vertex,label=below:$v_3$] (v3) at (2,-1) {};
					
					\draw[dedge] (v0)edge node[left,labelsty]{$e_1$}(v1);
					\draw[dedge] (v0)edge node[left,labelsty]{$e_2$}(v2);
					\draw[dedge] (v0)edge node[right,labelsty]{$e_3$}(v3);
				\end{tikzpicture}
				\caption{Non-loop to three neighbours.}
				\label{Fig1ExtensionNon3}
			\end{subfigure}
			\begin{subfigure}{0.49\textwidth}
				\centering
				\begin{tikzpicture}[scale=0.6]
					\node[vertex,label=below:$v_1$] (v1) at (-2,-1) {};
					\node[vertex,label=below:$v_2$] (v2) at (0,-1) {};
					
					\draw[dedge] (v1)edge node[above,labelsty]{$e$}(v2);
				\end{tikzpicture}
				{\Huge{$^{\mapsto}$}}
				\begin{tikzpicture}[scale=0.6]
					
					\node[vertex,label=above:$v_0$] (v0) at (0,1) {};
					\node[vertex,label=below:$v_1$] (v1) at (-2,-1) {};
					\node[vertex,label=below:$v_2$] (v2) at (2,-1) {};
					
					\draw[dedge] (v0)edge node[left,labelsty]{$e_1$}(v1);
					\draw[dedge] (v0)edge node[left,labelsty]{$e_2$}(v2);
					\draw[dedge] (v0)edge [bend left=40] node[right,labelsty]{$e_3$}(v2);
				\end{tikzpicture}
				\caption{Non-loop to two neighbours.}
				\label{Fig1ExtensionNon2}
			\end{subfigure}
            \begin{subfigure}{0.49\textwidth}
				\centering
				\begin{tikzpicture}[scale=1.2]
					\node[vertex,label=below:$v_1$] (v1) at (-2,-1) {};
					\node[vertex,label=below:$v_3$] (v3) at (0,-1) {};
					
					\draw[dedge] (v1)to [in=60,out=120,loop] node[above,labelsty]{$e$}(v1);
				\end{tikzpicture}
				{\Huge{$^{\mapsto}$}}
				\begin{tikzpicture}[scale=0.6]
					
					\node[vertex,label=above:$v_0$] (v0) at (0,1) {};
					\node[vertex,label=below:$v_1$] (v1) at (-2,-1) {};
					\node[vertex,label=below:$v_3$] (v3) at (2,-1) {};
					
					\draw[dedge] (v0)edge [bend right=40] node[left,labelsty]{$e_1$}(v1);
					\draw[dedge] (v0)edge node[right,labelsty]{$e_2$}(v1);
					\draw[dedge] (v0)edge node[right,labelsty]{$e_3$}(v3);
				\end{tikzpicture}
				\caption{Loop to two neighbours.}
				\label{Fig1ExtensionLoop2}
			\end{subfigure}
            \begin{subfigure}{0.49\textwidth}
				\centering
				\begin{tikzpicture}[scale=1.2]
					\node[vertex,label=below:$v_1$] (v1) at (-2,-1) {};
					\node[ivertex,label=below:\textcolor{white}{$v_2$}] (v2) at (-1,-1) {};

					\draw[dedge] (v1)to [in=60,out=120,loop] node[above,labelsty]{$e$}(v1);
				\end{tikzpicture}
				{\Huge{$^{\mapsto}$}}
				\begin{tikzpicture}[scale=0.6]
					
					\node[vertex,label=above:$v_0$] (v0) at (0,1) {};
					\node[vertex,label=below:$v_1$] (v1) at (0,-1) {};
					\draw[edge,white] (-1,0)edge node[left,labelsty]{\textcolor{black}{$e_1$}} (-1,0);
					\draw[edge,white] (1,0)edge node[right,labelsty]{\textcolor{black}{$e_3$}} (1,0);

                    \node[ivertex,label=below:\textcolor{white}{$v_1$}] (v1a) at (-2,-1) {};
					\node[ivertex,label=below:\textcolor{white}{$v_2$}] (v2a) at (2,-1) {};

					\draw[dedge] (v0)edge node[right,labelsty]{$e_2$}(v1);

					\draw[dedge] (v0) to [in=90,out=180](-1,0) to [in=180,out=270] (v1);
					\draw[dedge] (v0) to [in=90,out=0](1,0) to [in=0,out=270] (v1);

				\end{tikzpicture}
				\caption{Loop to one neighbour.}
				\label{Fig1ExtensionLoop1}
			\end{subfigure}
			\caption{The four variations of the gained $1$-extension.}
			\label{Fig1Extension}
		\end{figure}

		A \emph{gained loop-1-extension} forms a new gain graph $(G',m')$ by letting $m'|_{E(G)}=m$ and adding a new vertex $v_0$ with incident edges $l=(v_0,v_0;m'(l))$ and $e=(v_0,v_1;m'(e))$, for some $v_1\in V(G)$, with the additional requirement that $m'(l)$ has a non-trivial linear component. This move is illustrated in Figure \ref{FigLoop1Extension}. A \emph{gained loop-$1$-reduction} is the inverse move of a gained loop-$1$-extension.
		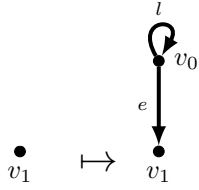
\begin{figure}[H]
			\centering
				\begin{tikzpicture}[scale=1.2]
					\node[vertex,white,label=right:\textcolor{white}{$v_0$}] (v0) at  (0,1) {};
					\node[vertex,label=below:$v_1$] (v1) at (0,0) {};
					\draw[dedge,white] (v0)to [in=60,out=120,loop] node[above,labelsty]{$l$}(v0);
					\draw[edge,white] (v0)to node[left,labelsty]{$e$}(v1);
				\end{tikzpicture}
				{\Huge{$^{\mapsto}$}}
				\begin{tikzpicture}[scale=1.2]
					\node[vertex,label=right:$v_0$] (v0) at  (0,1) {};
					\node[vertex,label=below:$v_1$] (v1) at (0,0) {};
					\draw[dedge] (v0)to [in=60,out=120,loop] node[above,labelsty]{$l$}(v0);
					\draw[dedge] (v0)to node[left,labelsty]{$e$}(v1);
				\end{tikzpicture}
			\caption{The gained loop-$1$-extension. Note that the gain $m'(l)$ must have a non-trivial linear component.}
			\label{FigLoop1Extension}
		\end{figure}
        
        A \emph{gained vertex-to-4-cycle move} forms a new gain graph $(G',m')$ by letting $m'|_{E(G)}=m$, choosing a vertex $v_1\in V(G)$ with incident edges $e_{12}=(v_1,v_2;m(e_{12}))$ and $e_{13}=(v_1,v_3;m(e_{13}))$ to distinct vertices $v_2,v_3\in V(G)$ and adding a new vertex $v_0$ with incident edges $e_{02}=(v_0,v_2;m'(e_{02}))$ and $e_{03}=(v_0,v_3;m'(e_{03}))$ such that $(m'(e_{02}))^{-1}m'(e_{03})=(m(e_{12}))^{-1}m(e_{13})$. Every other edge of the form $(v_1,w;m(e))\in E(G)$ (for some $w\in V(G)$) is either left unchanged or replaced by $(v_0,w;m(e))$. In the case where $w=v_1$, the edge may be left unchanged or be replaced by either $(v_0,v_1;m(e))$, $(v_1,v_0;m(e))$ or $(v_0,v_0;m(e))$. This move is illustrated in Figure \ref{Fig4Cyclextension}. A \emph{gained $4$-cycle-to-vertex move} is the inverse move of a gained vertex-to-$4$-cycle move.
		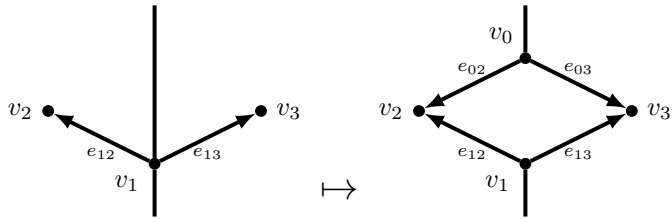
\begin{figure}[H]
			\centering
			\begin{tikzpicture}[scale=0.7]
				\node[vertex,label=200:$v_1$] (v1) at (0,-1) {};
				\node[vertex,label=left:$v_2$] (v2) at (-2,0) {};
				\node[vertex,label=right:$v_3$] (v3) at (2,0) {};
				
				\draw[dedge] (v1)edge node[below,labelsty]{$e_{12}$}(v2);
				\draw[dedge] (v1)edge node[below,labelsty]{$e_{13}$}(v3);
				\draw[edge] (v1)edge(0,-2);
				\draw[edge] (v1)edge(0,2);
			\end{tikzpicture}
			{\Huge{$^{\mapsto}$}}
			\begin{tikzpicture}[scale=0.7]
				\node[vertex,label=110:$v_0$] (v0) at (0,1) {};
				\node[vertex,label=200:$v_1$] (v1) at (0,-1) {};
				\node[vertex,label=left:$v_2$] (v2) at (-2,0) {};
				\node[vertex,label=right:$v_3$] (v3) at (2,0) {};
				
				\draw[dedge] (v1)edge node[below,labelsty]{$e_{12}$}(v2);
				\draw[dedge] (v1)edge node[below,labelsty]{$e_{13}$}(v3);
				\draw[dedge] (v0)edge node[above,labelsty]{$e_{02}$}(v2);
				\draw[dedge] (v0)edge node[above,labelsty]{$e_{03}$}(v3);
				\draw[edge] (v1)edge(0,-2);
				\draw[edge] (v0)edge(0,2);
			\end{tikzpicture}
			\caption{The gained vertex-to-$4$-cycle move.}
			\label{Fig4Cyclextension}
		\end{figure}
        
        A \emph{gained vertex-to-$K_4$ move} forms a new gain graph $(G',m')$ by letting $m'|_{E(G)}=m$ and replacing a vertex $v_1\in V(G)$ with a copy of $K_4$ in which all edges have identity gain. Every other edge of the form $(v_1,w;m(e_w))\in E(G)$ (for some $w\in V(G)$ which may be equal to $v_1$) is either left unchanged or replaced by $(v_0,w;m(e_w))$, where $v_0$ is another vertex in the copy of $K_4$. Note that different edges of this form may use different choices of $v_0$. In the case where $w=v_1$, either end vertex can be replaced by any vertex in the copy of $K_4$.  This move is illustrated in Figure \ref{FigK4extension}. A \emph{gained $K_4$-to-vertex move} is the inverse move of a gained vertex-to-$K_4$ move. 
		
		\begin{figure}[H]
			\centering
			\begin{tikzpicture}[scale=0.7]
				\node[vertex,label=left:$v_1$] (v1) at (0,0) {};
				
				\draw[edge] (v1)edge(-1,-1);
				\draw[edge] (v1)edge(1,-1);
			\end{tikzpicture}
			{\Huge{$^{\mapsto}$}}
			\begin{tikzpicture}[scale=0.7]
				\node[vertex,label=left:$v_1$] (v1) at (0,0) {};
				\node[vertex] (w2) at (1,0) {};
				\node[vertex] (w3) at (1,1) {};
				\node[vertex] (w4) at (0,1) {};
				
				\draw[edge] (v1)edge(w2);
				\draw[edge] (w2)edge(w3);
				\draw[edge] (w3)edge(w4);
				\draw[edge] (w4)edge(v1);
				\draw[edge] (v1)edge(w3);
				\draw[edge] (w4)edge(w2);
				
				\draw[edge] (v1)edge(-1,-1);
				\draw[edge] (w2)edge(2,-1);
			\end{tikzpicture}
			\caption{The gained vertex-to-$K_4$ move.}
			\label{FigK4extension}
		\end{figure}
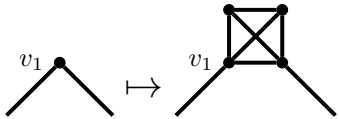

        A \emph{gained edge-to-$K_3$ move} (also known as a \emph{vertex-splitting move}) forms a new gain graph $(G',m')$ by letting $m'|_{E(G)}=m$, choosing an edge $e=(v_1,v_2;m(e))\in E(G)$ and adding a new vertex $v_0$ with incident edges $e_1=(v_0,v_1;m'(e_1))$ and $e_2=(v_0,v_2;m'(e_2))$ such that $(m'(e_1))^{-1}m'(e_2)=m(e)$. Every other edge of the form $(v_1,w;m(e))\in E(G)$ (for some $w\in V(G)$) is either left unchanged or replaced by $(v_0,w;m(e))$. In the case where $w=v_1$, the edge may be left unchanged or be replaced by either $(v_0,v_1;m(e))$, $(v_1,v_0;m(e))$ or $(v_0,v_0;m(e))$. This move is illustrated in Figure \ref{FigVertexSplitting}. A \emph{gained $K_3$-to-edge move} is the inverse move of a gained edge-to-$K_3$ move.
		\begin{figure}[H]
			\centering
			\begin{tikzpicture}[scale=0.7]
				\node[vertex,label=below:$v_1$] (v1) at (0,0) {};
				\node[vertex,label=below:$v_2$] (v2) at (3,0) {};
				
				\draw[edge] (v1)edge node[below,labelsty]{$e$}(v2);
				\draw[edge] (v1)edge(-1,-1);
				\draw[edge] (v1)edge(-1,1);
			\end{tikzpicture}
			{\Huge{$^{\mapsto}$}}
			\begin{tikzpicture}[scale=0.7]
				\node[vertex,label=above:$v_0$] (v0) at (0,2) {};
				\node[vertex,label=below:$v_1$] (v1) at (0,0) {};
				\node[vertex,label=below:$v_2$] (v2) at (3,0) {};
				
				\draw[dedge] (v1)edge node[below,labelsty]{$e$}(v2);
				\draw[dedge] (v0)edge node[left,labelsty]{$e_1$}(v1);
				\draw[dedge] (v0)edge node[above,labelsty]{$e_2$}(v2);
				\draw[edge] (v1)edge(-1,-1);
				\draw[edge] (v0)edge(-1,3);
			\end{tikzpicture}
			\caption{The gained edge-to-$K_3$ move.}
			\label{FigVertexSplitting}
		\end{figure}
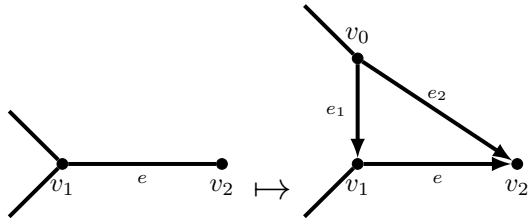

        Let $(G_1,m_1)$ and $(G_2,m_2)$ be disjoint $\Gamma$-gain graphs. An \emph{edge-joining move} forms a new $\Gamma$-gain graph $(G_1\oplus G_2,m')$ by taking the disjoint union $G_1\cup G_2$, letting $m'|_{E(G_1)}=m_1$, letting $m'|_{E(G_2)}=m_2$ and adding an edge $e=(v_1,v_2;m'(e))$, where $v_1\in V(G_1)$ and $v_2\in V(G_2)$. This creates the $\Gamma$-gain graph
		\begin{align*}
			G_1\oplus G_2=(V(G_1)\cup V(G_2),E(G_1)\cup E(G_2)\cup((v_1,v_2;m'(e)))).
		\end{align*}
        This move is illustrated in Figure \ref{FigEdgeJoining}. 
        \begin{figure}[H]
			\centering
			\begin{tikzpicture}[scale=0.2]
				\node[vertex,label=below:$v_1$] (x) at (-4,0) {};
				\node[vertex,label=below:$v_2$] (y) at (4,0) {};

				\draw[dashed](-7,0)circle(5);
				\draw[dashed](7,0)circle(5);
                \node[blankvertex] (G1) at (-10,0) {$G_1$};
                \node[blankvertex] (G2) at (10,0) {$G_2$};
			\end{tikzpicture}
			{\Huge{$^{\mapsto}$}}
			\begin{tikzpicture}[scale=0.2]
				\node[vertex,label=below:$v_1$] (x) at (-4,0) {};
				\node[vertex,label=below:$v_2$] (y) at (4,0) {};

				\draw[edge] (x)  to  (y);

				\draw[dashed](-7,0)circle(5);
				\draw[dashed](7,0)circle(5);
                \node[blankvertex] (G1) at (-10,0) {$G_1$};
                \node[blankvertex] (G2) at (10,0) {$G_2$};
			\end{tikzpicture}
			\caption{The edge-joining move}
			\label{FigEdgeJoining}
		\end{figure}
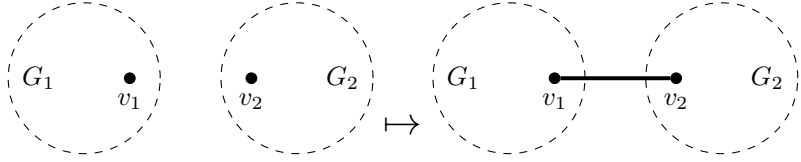
	\end{definition}
   
    \section{Rigidity in $\ell_q$-planes for $q\in(1,\infty)\backslash\{2\}$}
    \subsection{Preliminaries on Rigidity and Symmetric Rigidity in $\ell_q$-planes for $q\in(1,\infty)\backslash\{2\}$}\label{SectionLqBasics}
    Given some $q\in[1,\infty)$, the \emph{$\ell_q$-norm} $\|\cdot\|_q$ on $\R^2$ is defined, for all $a=(a_1,a_2)\in\R^2$, by
    \begin{align*}
        \|a\|_q=\sqrt[q]{|a_1|^q+|a_2|^q}.
    \end{align*}
    Note that the $\ell_2$-norm is simply the familiar Euclidean norm on $\R^2$.
    
    In this section, we briefly summarise the basics of rigidity in the space $(\R^2,\|\cdot\|_q)$ for $q\in(1,\infty)\backslash\{2\}$. For more details, see \cite{NonEuclidean}. Forced-symmetric and periodic rigidity in this setting will be studied in Sections \ref{SectionLqPeriodic} to \ref{SectionLqWallpaper}. Note that studying rigidity in the $\ell_1$-plane or the $\ell_\infty$-plane requires different techniques, which will be discussed in Sections \ref{SectionPolytopicBasics} to \ref{SectionPolytopicWallpaper}.
    
    For a point $a=(a_1,a_2)\in\R^2$ and a scalar $k\in(0,\infty)$, define
    \begin{align*}
        a^{(k)}=(\mathrm{sgn}(a_1)|a_1|^k,\mathrm{sgn}(a_2)|a_2|^k).
    \end{align*}
     We can now see how the general definition of an infinitesimal motion applies with these norms, as proved in \cite[Proposition 3.2]{NonEuclidean}. Let $G$ be a simple undirected graph and let $(G,p)$ be a framework in $(\R^2,\|\cdot\|_q)$ for some $q\in(1,\infty)\backslash\{2\}$. An \emph{infinitesimal motion} of $(G,p)$ is a map $u:V(G)\to\R^2$ such that, for all edges $\{v_i,v_j\}\in E(G)$,
    \begin{align*}
        (p(v_i)-p(v_j))^{(q-1)}\cdot(u(v_i)-u(v_j))=0.
    \end{align*}
    We can also define the rigidity matrix for a framework in these spaces, which is analogous to the standard rigidity matrix for Euclidean spaces. The \emph{rigidity matrix} $R_q(G,p)$ of $(G,p)$ is the $|E(G)|\times 2|V(G)|$ matrix with one row for each edge in $E(G)$ and a $2$-tuple of columns for each vertex in $V(G)$, where the row for an edge $e=\{v_i,v_j\}\in E(G)$ is
     \setcounter{MaxMatrixCols}{20}
		\begin{align*}
        \begin{bNiceMatrix}[first-row,first-col]
         &&&&v_i&&&&v_j&&&&\\
		  e&0&...&0&(p(v_i)-p(v_j))^{(q-1)}&0&...&0&(p(v_j)-p(v_i))^{(q-1)}&0&...&0
        \end{bNiceMatrix}
        .
    \end{align*}
    From this, it is clear that $u\in\R^{2|V(G)|}$ is an infinitesimal motion for $(G,p)$ if and only if $R_q(G,p)u=0$, so the kernel of the rigidity matrix represents the space of infinitesimal motions \cite[Proposition 3.2]{NonEuclidean}.

    We now define the orbit rigidity matrix for the general $\ell_q$-plane. This is analogous to the Euclidean orbit rigidity matrix, which is described in \cite{Orbit} for finite symmetries and \cite{Torus} for translation symmetries. Let $\Gamma$ be a group of isometries on $(\R^2,\|\cdot\|_q)$ and let $(G,m)$ be a $\Gamma$-gain graph with a configuration $p:V(G)\to\R^2$ that derives a $\Gamma$-symmetric framework $(\G,\p)$ in $(\R^2,\|\cdot\|_q)$. The \emph{orbit rigidity matrix} $\mathbf{O}_q(\G,\p,\Gamma)$ of $(\G,\p)$ is the $|E(G)|\times 2|V(G)|$ matrix with one row for each edge in $(G,m)$ and a $2$-tuple of columns for each vertex in $(G,m)$, defined as follows. There are two different possibilities for the row corresponding to an edge $e=(v_i,v_j;m(e))\in E(G)$:
		\begin{enumerate}
			\item Suppose that $e$ is not a loop $(v_i\neq v_j)$. Then the row for $e$ is
            \setcounter{MaxMatrixCols}{20}
		\begin{align*}
        \begin{bNiceMatrix}[first-row,first-col]
         &&&&v_i&&&&v_j&&&&\\
		  e&0&...&0&(p(v_i)-m(e)p(v_j))^{(q-1)}&0&...&0&(p(v_j)-(m(e))^{-1}p(v_i))^{(q-1)}&0&...&0
        \end{bNiceMatrix}
        .
    \end{align*}
			\item Suppose that $e$ is a loop $(v_i=v_j)$. Then the row for $e$ is
            \begin{align*}
        \begin{bNiceMatrix}[first-row,first-col]
         &&&&v_i&&&&\\
         e&0&...&0&(p(v_i)-m(e)p(v_i))^{(q-1)}+(p(v_i)-(m(e))^{-1}p(v_i))^{(q-1)}&0&...&0
        \end{bNiceMatrix}
        .
    \end{align*}
		\end{enumerate}
    As in the Euclidean plane, the kernel of the orbit rigidity matrix corresponds to the space of symmetric infinitesimal motions. The proof is an immediate extension of that seen for the Euclidean orbit rigidity matrix in \cite{Orbit}. See \cite{EssonThesis} for details.

    \subsection{Periodic Rigidity in $\ell_q$-planes for $q\in(1,\infty)\backslash\{2\}$}\label{SectionLqPeriodic}

    In this section, we fill an important gap in the theory, as we study forced-periodic rigidity in non-Euclidean $\ell_q$-planes for the first time, although the proof uses a considerable amount of previous work. Note that we are only considering fixed-lattice periodic rigidity here. When considering periodic frameworks in any plane, we can assume, without loss of generality, that the periodicity lattice vectors are $(1,0)$ and $(0,1)$. Put otherwise, we are considering $\Gamma$-symmetric frameworks, where $\Gamma$ is the group of isometries generated by translations by $(1,0)$ and $(0,1)$. For ease, we identify this group with $\Z^2$. The following result characterises periodic rigidity in the  $\ell_q$-plane for $q\in(1,\infty)\backslash\{2\}$. 
    \begin{theorem}\label{lqPeriodic}
        Let $q\in(1,\infty)\backslash\{2\}$. Let $(\G,\p)$ be a $\Z^2$-regular framework in $(\R^2,\|\cdot\|_q)$ with underlying $\Z^2$-gain graph $(G,m)$. Then $(\G,\p)$ is minimally periodically infinitesimally rigid if and only if $(G,m)$ is $(2,2)$-tight.
    \end{theorem}
    
    The necessity part of Theorem \ref{lqPeriodic} is proved in the standard way, see e.g. \cite{NonEuclidean}, noting that the $\ell_q$-plane admits a $2$-dimensional space of trivial periodic infinitesimal motions (induced by all translations) and therefore the required rank is $2|V(G)|-2$. It can then be shown that any subgraph $G'\subseteq G$ with $|E(G')|>2|V(G')|-2$ creates a row-dependence in the orbit rigidity matrix, contradicting minimal rigidity.

    To prove the sufficiency part of Theorem \ref{lqPeriodic}, we use an inductive approach, which is similar to the proofs of \cite[Theorem 3.6]{NonEuclidean} and \cite[Theorem 5.1]{Inductive}. The inductive method involves gained $0$-extensions, $1$-extensions, vertex-to-$4$-cycle moves and vertex-to-$K_4$ moves. The first step is to show that each of these extensions preserves minimal rigidity of $\Z^2$-gain graphs in $(\R^2,\|\cdot\|_q)$. This can be done using straightforward adaptations of standard  methods.
    \begin{proposition}\label{PeriodicLqExtensions}
		Let $(\G,\p)$ be a periodic framework in $(\R^2,\|\cdot\|_q)$ that is minimally periodically infinitesimally rigid. Let $(G,m)$ be the underlying $\Z^2$-gain graph of $(\G,\p)$. Let $(G',m')$ be formed by a gained $0$-extension, $1$-extension, vertex-to-$4$-cycle move or vertex-to-$K_4$ move of $(G,m)$. Let $(\G',\p')$ be a $\Z^2$-regular framework derived from $(G',m')$. Then $(\G',\p')$ is minimally periodically infinitesimally rigid.
	\end{proposition}
    \begin{proof}
        The fact that a gained $0$-extension preserves minimal periodic rigidity can be proved by a slight modification of the method seen in \cite[Lemma 3.8]{NonEuclidean}. This involves choosing a configuration in which the new vertex is not collinear with its derived neighbours. Then basic linear algebra can be used to show that the orbit rigidity matrix is row-independent.

        The fact that a gained $1$-extension preserves minimal periodic rigidity can be proved by combining the standard methods seen in \cite[Lemma 3.10]{NonEuclidean} and \cite[Proposition 3.2]{Inductive}. After choosing a configuration in which the new vertex is placed on the line corresponding to the edge that was removed, basic linear algebra can be used to show that the orbit rigidity matrix is row-independent. Note that $(2,2)$-tight multigraphs have no loops, so the $1$-extension variations seen in Figures \ref{Fig1ExtensionLoop2} and \ref{Fig1ExtensionLoop1} are not relevant here.

        Vertex-to-$4$-cycle moves and vertex-to-$K_4$ moves can be proved to preserve minimal periodic rigidity by adapting the methods used to prove \cite[Lemma 3.12]{NonEuclidean} and \cite[Lemma 3.14]{NonEuclidean} respectively. We refer the reader to \cite{EssonThesis} for details.
    \end{proof}

    The following result provides the inductive construction of $(2,2)$-tight $\Z^2$-gain graphs.
    \begin{theorem}\cite[Theorem 15]{Surfaces}\label{PeriodicLqInduction}
		A $\Z^2$-gain graph is $(2,2)$-tight if and only if it can be formed from $K_1$ by a sequence of gained $0$-extensions, $1$-extensions, vertex-to-$4$-cycle moves and vertex-to-$K_4$ moves.
	\end{theorem}
    This allows the proof of sufficiency for Theorem \ref{lqPeriodic} to be completed.
    \begin{proof}[Proof of sufficiency for Theorem \ref{lqPeriodic}]
        It is clear that any framework derived from a $\Z^2$-gain graph on $K_1$ is minimally periodically infinitesimally rigid in this space.  By Theorem \ref{PeriodicLqInduction}, every $(2,2)$-tight $\Z^2$-gain graph can be formed from $K_1$ by a sequence of gained $0$-extensions, $1$-extensions, vertex-to-$4$-cycle moves and vertex-to-$K_4$ moves. Proposition \ref{PeriodicLqExtensions} shows that all of these moves preserve minimal periodic infinitesimal rigidity of $\Z^2$-regular derived frameworks. Hence, every $\Z^2$-regular framework in $(\R^2,\|\cdot\|_q)$ derived from a $(2,2)$-tight $\Z^2$-gain graph is minimally periodically infinitesimally rigid.
    \end{proof}

    \subsection{Reflective Symmetric Rigidity in $\ell_q$-planes for $q\in(1,\infty)\backslash\{2\}$}\label{SectionLqReflection}
    In this section, we present a combinatorial characterisation for minimal reflectionally-symmetric infinitesimal rigidity in the $\ell_q$-plane for $q\in(1,\infty)\backslash\{2\}$. Let $\mathcal{C}_s=\{0,s\}$ be a reflection symmetry group. Note that, in these planes, symmetry groups must be compatible with the symmetries of the unit ball. Hence, the linear part of the symmetry group must be a subset of the dihedral group of order $8$. This means that axes of reflection must be either the $x$-axis, the $y$-axis or one of the diagonals. Since the geometric results here do not rely on the position of the reflection line, we can, without loss of generality, assume that the reflection $s$ is in the $x$-axis.

    \begin{theorem}\label{Reflectivelq}
        Let $q\in(1,\infty)\backslash\{2\}$. Let $(\G,\p)$ be a $\mathcal{C}_s$-regular framework in $(\R^2,\|\cdot\|_q)$ with underlying $\mathcal{C}_s$-gain graph $(G,m)$. Then $(\G,\p)$ is minimally $\mathcal{C}_s$-symmetrically infinitesimally rigid if and only if $(G,m)$ is $(2,2,1)$-gain-tight.
    \end{theorem}
	Necessity of this condition is again proved in the standard way using the new orbit rigidity matrix. Note that there is a $1$-dimensional space of trivial $\mathcal{C}_s$-symmetric infinitesimal motions (induced by translations parallel to the axis of reflection), so the required overall count for the gain graph is $|E(G)|=2|V(G)|-1$. It can be seen that any subgraph that breaks the sparsity conditions will create a row-dependence in the orbit rigidity matrix.
    
    We now prove sufficiency. Again, we prove this inductively. The relevant extension types are gained $0$-extensions, $1$-extensions, loop-$1$-extensions, vertex-to-$4$-cycle moves, vertex-to-$K_4$ moves and edge-joining moves. The first four of these can easily be proved to preserve minimal $\mathcal{C}_s$-symmetric rigidity of $\mathcal{C}_s$-regular derived frameworks by adapting standard methods, such as those seen in \cite{EGRES} and \cite{NonEuclidean}.
    \begin{proposition}\label{LqCsStraightforwardExtensions}
		Let $(\G,\p)$ be a $\mathcal{C}_s$-symmetric framework in $(\R^2,\|\cdot\|_q)$ that is minimally $\mathcal{C}_s$-symmetrically infinitesimally rigid. Let $(G,m)$ be the underlying $\mathcal{C}_s$-gain graph of $(\G,\p)$. Let $(G',m')$ be formed by a gained $0$-extension, $1$-extension, loop-$1$-extension or vertex-to-$4$-cycle move of $(G,m)$. Let $(\G',\p')$ be a $\mathcal{C}_s$-regular framework derived from $(G',m')$. Then $(\G',\p')$ is minimally $\mathcal{C}_s$-symmetrically infinitesimally rigid.
	\end{proposition}
    To prove that vertex-to-$K_4$ moves preserve minimal rigidity, we adapt a method from \cite[Lemma 9]{Surfaces}. Note that the method used in \cite[Lemma 3.14]{NonEuclidean} does not work for $\mathcal{C}_s$-symmetric frameworks, as vertical translations are not $\mathcal{C}_s$-symmetric.
    \begin{proposition}\label{LqCsK4Extension}
		Let $(\G,\p)$ be a $\mathcal{C}_s$-symmetric framework in $(\R^2,\|\cdot\|_q)$ that is minimally $\mathcal{C}_s$-symmetrically infinitesimally rigid. Let $(G,m)$ be the underlying $\mathcal{C}_s$-gain graph of $(\G,\p)$. Let $(G',m')$ be formed by a gained vertex-to-$K_4$ move of $(G,m)$. Let $(\G',\p')$ be a $\mathcal{C}_s$-regular framework derived from $(G',m')$. Then $(\G',\p')$ is minimally $\mathcal{C}_s$-symmetrically infinitesimally rigid.
	\end{proposition}
    \begin{proof}
        Note first that the move preserves $(2,1)$-tightness, so it is sufficient to show that the move preserves rigidity. Suppose that the vertex-to-$K_4$ move replaces $v_1\in V(G)$ with the balanced subgraph $K_4(v_1,w_2,w_3,w_4)$. By Lemma \ref{lemma:swicth}, it can be assumed that every edge in this copy of $K_4$ has identity gain. For contradiction, suppose that $(\G',\p')$ is $\mathcal{C}_s$-symmetrically infinitesimally flexible.
		
		Let $p:V(G)\to\R^{2|V(G)|}$ be a $\mathcal{C}_s$-regular configuration of $(G,m)$. Define a sequence of $\mathcal{C}_s$-regular configurations $(p^{k})_{k\in\N}$ of $(G',m')$ with $p^{k}|_{V(\G)}=p$ for all $k\in\N$. Let $\p^{k}$ be the derived configuration of each $p^{k}$. Suppose that, for $i\in\{2,3,4\}$, we have that $p^{k}(w_i)\to p(v_1)$ as $k\to\infty$. Hence, $(p^{k})_{k\in\N}$ is convergent. Let $p^\infty$ denote its limit, deriving a configuration $\p^\infty$ of $\G'$.
		
		Since $(\G',\p)$ is assumed to be $\mathcal{C}_s$-regular and $\mathcal{C}_s$-symmetrically infinitesimally flexible, it must be that each $(\G',\p^{k})$ is also $\mathcal{C}_s$-symmetrically infinitesimally flexible. For each $k\in\N$, let $u^k\in\R^{2|V(G')|}$ be a unit (with respect to the Euclidean norm) vector that derives a non-trivial $\mathcal{C}_s$-symmetric infinitesimal motion of $(\G',\p^{k})$. This can be chosen so that it is orthogonal (in the Euclidean sense) to the $1$-dimensional space of trivial infinitesimal motions. Since the sequence $(u^k)_{k\in\N}$ is formed only of unit vectors, this is a bounded sequence and it therefore has a convergent subsequence by the Bolzano-Weierstrass Theorem. Let $u\in\R^{2|V(G')|}$ be the limit of this convergent subsequence. We aim to show that $u$ derives an infinitesimal motion of $(\G',\p^\infty)$. By passing to a subsequence if necessary, it can be assumed that $u^k\to u$.
        
        As seen in \cite[Example 3.4]{NonEuclidean}, any $\mathcal{C}_s$-regular framework derived from $K_4$ is infinitesimally rigid. The restriction of each $u^k$ to the copy of $K_4$ is therefore a trivial infinitesimal motion. This means that, for all $k\in\N$, $u^k(v_1)=u^k(w_2)=u^k(w_3)=u^k(w_4)$. Hence, $u(v_1)=u(w_2)=u(w_3)=u(w_4)$.
		
		Since the sequence $(p^{k}(v))_{k\in\N}$ is constant for all $v\in V(\G)$, the infinitesimal motion constraints given by each edge between vertices of $V(G)$ are unchanged as $k\to\infty$. Bars of zero length do not impose any constraints on infinitesimal motions, so the edges that are part of $K_4(v_1,w_2,w_3,w_4)$ need not be considered for $u$.
        
        Consider an edge $e=(w,v;m'(e))\in E(G')$, where $w\in\{w_2,w_3,w_4\}$ and $v\in V(G)\backslash\{v_1\}$. The constraint placed on $u^k$ by this edge is
		\begin{align*}
			(p^{k}(w)-m'(e)p(v))^{(q-1)}\cdot(u^k(w)-m'(e)u^k(v))=0.
		\end{align*}
		Since $p^{k}(w)\to p(v_1)$ and $u^k(w)\to u(v_1)$ and $u^k(v)\to u(v)$, it can be seen that
		\begin{align*}
			(p^{k}(w)-m'(e)p(v))^{(q-1)}\cdot(u^k(w)-m'(e)u^k(v))\to(p(v_1)-m'(e)p(v))^{(q-1)}\cdot(u(v_1)-m'(e)u(v)).
		\end{align*}
		Since the sequence on the left is constrained to be a constant zero sequence, it can be seen that
		\begin{align*}
			(p(v_1)-m'(e)p(v))^{(q-1)}\cdot(u(v_1)-m'(e)u(v))=0.
		\end{align*}
		Hence, $u$ satisfies the infinitesimal motion constraint imposed by $e$.
		
		It remains to consider edges of non-identity gain between vertices within the copy of $K_4$, which arise from loops at $v_1$ under the extension. Let $e=(w_i,w_j;m'(e))\in E(G')$ be such an edge for some distinct pair $w_i,w_j\in\{v_1,w_2,w_3,w_4\}$. This imposes the following constraint on $u^k$:
		\begin{align*}
			(p^{k}(w_i)-m'(e)p^{k}(w_j))^{(q-1)}\cdot(u^k(w_i)-m'(e)u^k(w_j))=0.
		\end{align*}
		Like before, note that
		\begin{align*}
			(p^{k}(w_i)-m'(e)p^{k}(w_j))^{(q-1)}\cdot(u^k(w_i)-m'(e)u^k(w_j))\to(p(v_1)-m'(e)p(v_1))\cdot(u(v_1)-m'(e)u(v_1)).
		\end{align*}
		Since the sequence on the left is constrained to be a constant zero sequence,
		\begin{align*}
			(p(v_1)-m'(e)p(v_1))\cdot(u(v_1)-m'(e)u(v_1))=0.
		\end{align*}
		Hence, $u$ satisfies the infinitesimal motion constraint imposed by $e$. Hence, $u$ is an infinitesimal motion of $(\G',\p^\infty)$.
		
		Note that each $u^k$ is a non-trivial infinitesimal motion that has unit norm and is orthogonal to every trivial infinitesimal motion. For any trivial infinitesimal motion $x$ and any $k\in\N$, we have that $u^k\cdot x=0$. By continuity of the dot product, $u\cdot x=0$. By continuity of the norm, $\|u\|_2=1$. Hence, $u$ derives a non-trivial $\mathcal{C}_s$-symmetric infinitesimal motion of $(\G',\p^\infty)$. Since $u$ is constant on $\{v_1,w_2,w_3,w_4\}$, it can be seen that the restriction $u|_{V(G)}$ derives a non-trivial $\mathcal{C}_s$-symmetric infinitesimal motion on the $\mathcal{C}_s$-regular framework $(\G,\p)$, which contradicts the fact that this is $\mathcal{C}_s$-symmetrically infinitesimally rigid. Hence, vertex-to-$K_4$ moves preserve minimal $\mathcal{C}_s$-symmetric infinitesimal rigidity.
    \end{proof}
    We now show that edge-joining moves preserve minimal $\mathcal{C}_s$-symmetric infinitesimal rigidity, following the idea seen in \cite[Lemma 11]{Surfaces}.
    \begin{proposition}\label{LqCsJoiningExtension}
		Let $(\G_1,\p_1)$ and $(\G_2,\p_2)$ be $\mathcal{C}_s$-symmetric frameworks in $(\R^2,\|\cdot\|_q)$ that are minimally $\mathcal{C}_s$-symmetrically infinitesimally rigid. Let $(G_1,m_1)$ and $(G_2,m_2)$ be the underlying $\mathcal{C}_s$-gain graphs of $(\G_1,\p_1)$ and $(\G_2,\p_2)$ respectively. Let $(G',m')$ be formed by a gained edge-joining move that combines $(G_1,m')$ and $(G_2,m')$. Let $(\G',\p')$ be a $\mathcal{C}_s$-regular framework derived from $(G',m')$. Then $(\G',\p')$ is minimally $\mathcal{C}_s$-symmetrically infinitesimally rigid.
	\end{proposition}
	\begin{proof}
		Suppose that the edge-joining move gives $G'= G_1\oplus G_2=(V(G_1)\cup V(G_2),E(G_1)\cup E(G_2)\cup\{(v_1,v_2;m'(e))\})$, where $v_1\in V(G_1)$ and $v_2\in V(G_2)$. It is clear that this move preserves $(2,1)$-tightness, so it remains only to prove that it preserves row-independence of the orbit rigidity matrix. Since translating a framework does not change its rigidity properties, without loss of generality, we may assume that $\tilde{p}_1(v_1)\neq m'(e)\tilde{p}_2(v_2)$, and  we may further assume that $\p'|_{V(\G_1)}=\p_1$ and $\p'|_{V(\G_2)}=\p_2$. If we prove that a framework following these assumptions is $\mathcal{C}_s$-symmetrically infinitesimally rigid, then $\mathcal{C}_s$-symmetric infinitesimal rigidity follows for all $\mathcal{C}_s$-regular configurations of $\G$. Observe that the orbit rigidity matrix of $(\G',\p')$ takes the following form:
		\begin{align*}
			\mathbf{O}_q(\G',\p',\mathcal{C}_s)=
			\begin{bmatrix}
				\mathbf{O}_q(\G_1,\p'|_{V(\G_1)},\mathcal{C}_s)&0\\
				0&\mathbf{O}_q(\G_2,\p'|_{V(\G_2)},\mathcal{C}_s)\\
				\begin{matrix}
					(\p'(v_1)-m'(e)\p'(v_2))^{(q-1)}&0
				\end{matrix}
				&
				\begin{matrix}
					(\p'(v_2)-(m'(e))^{-1}\p'(v_1))^{(q-1)}&0
				\end{matrix}
			\end{bmatrix}
			.
		\end{align*}
		It is clear that the rows for $E(G_1)\cup E(G_2)$ are linearly independent, so it remains only to show that the row for the new edge $e$ is independent of the others.

        The rank of $\mathbf{O}_q(\G_1,\p'|_{V(\G_1)},\mathcal{C}_s)\oplus\mathbf{O}_q(\G_2,\p'|_{V(\G_2)},\mathcal{C}_s)$ is $2|V(G)|-2$, so the resulting kernel is $2$-dimensional. One element of this kernel is the vector $u\in \R^{2|V(G)|}$ defined by setting $u(v)=(1,0)$ if $v\in V(G_1)$ and $u(v)=(0,0)$ if $v\in V(G_2)$. By regularity, it can be assumed that $\p'(v_1)-m'(e)\p'(v_2)$ has a non-zero horizontal component. From this, it follows that $e$ eliminates $u$ from the kernel, which leaves only a $1$-dimensional kernel, as required.
	\end{proof}
    Having seen that each of these extensions preserves minimal rigidity, consider the following theorem.
	\begin{theorem}\label{Induction221}\cite[Theorem 17]{Surfaces}
		A $\mathcal{C}_s$-gain graph is $(2,2,1)$-gain-tight if and only if it can be constructed from either an unbalanced gain graph on $K_1^1$ (the multigraph consisting of a single vertex with a single loop)  or a gain graph on $K_4+e$, where every edge except $e$ has identity gain, by a sequence of gained $0$-extensions, $1$-extensions, loop-$1$-extensions, vertex-to-$K_4$ moves, vertex-to-$4$-cycle moves and edge-joining moves.
	\end{theorem}
    This allows the proof of sufficiency for Theorem \ref{Reflectivelq} to be completed.
	\begin{proof}[Proof of sufficiency for Theorem \ref{Reflectivelq}]
		It is easy to see that all $\mathcal{C}_s$-regular frameworks derived from gain graphs on $K_1^1$ are minimally $\mathcal{C}_s$-symmetrically infinitesimally rigid. The same applies for $\mathcal{C}_s$-regular frameworks derived from $\mathcal{C}_s$-gain graphs on $K_4+e$, where every edge except $e$ has identity gain. Theorem \ref{Induction221} states that every $(2,2,1)$-gain-tight graph can be constructed from one of these base graphs by a sequence of gained $0$-extensions, $1$-extensions, loop-$1$-extensions, vertex-to-$K_4$ moves, vertex-to-$4$-cycle moves and edge-joining moves. Propositions \ref{LqCsStraightforwardExtensions}, \ref{LqCsK4Extension} and \ref{LqCsJoiningExtension} show that all of these moves preserve minimal $\mathcal{C}_s$-symmetric infinitesimal rigidity of $\mathcal{C}_s$-regular derived frameworks. Hence, every $\mathcal{C}_s$-regular framework in $(\R^2,\|\cdot\|_q)$ derived from a $(2,2,1)$-gain-tight $\mathcal{C}_s$-gain graph is minimally $\mathcal{C}_s$-symmetrically infinitesimally rigid.
	\end{proof}

    \subsection{$\Z^2\rtimes\mathcal{C}_s$-Symmetric Rigidity in $\ell_q$-planes for $q\in(1,\infty)\backslash\{2\}$}\label{SectionLqWallpaper}
    We now consider frameworks that are symmetric with respect to the wallpaper group $\Z^2\rtimes\mathcal{C}_s$. This is the wallpaper group formed by taking the semi-direct product of the group $\Z^2$ of translations (w.l.o.g. generated by the vectors $(1,0)$ and $(0,1)$) with the reflectional group $\mathcal{C}_s$, generated by the reflection $s$, in the plane (where w.l.o.g. the mirror line of $s$ is the $x$-axis) \cite[Section 3.2]{Gocg}. In the Hermann–Mauguin notation used in crystallography, this group is  denoted by $pm$ \cite{Schatt, HMbook}, while in the orbifold notation advocated by J.H. Conway, this group is known as $\ast\ast$ \cite{Sym}.
    
    For ease, we will use the following notation for elements of $\Z^2\rtimes\mathcal{C}_s$. Let $a,b\in\Z$ and let $r\in\mathcal{C}_s$. Then $(a,b,r)\in\Z^2\rtimes\mathcal{C}_s$ is the isometry formed by composing the reflection $r$, followed by a translation by the vector $(a,b)$.
    
    In \cite[Theorem 3.2]{EssonCrystallographic}, we characterised conditions for $\Z^2\rtimes\mathcal{C}_s$-symmetric rigidity in the Euclidean plane. In this section, we do the same in non-Euclidean $\ell_q$-planes with $q\in(1,\infty)\backslash\{2\}$. As before, we are assuming that the periodicity lattice must remain fixed.
    \begin{definition}\label{DefZReflectivelq}
        A $\Z^2\rtimes\mathcal{C}_s$-gain graph $(G,m)$ is said to be \emph{$(\Z^2\rtimes\mathcal{C}_s)_q$-tight} if it satisfies the following conditions:
        \begin{enumerate}
			\item \emph{$(2,1)$-tight Condition:} $G$ is $(2,1)$-tight;
			\item \emph{Purely Periodic Condition:} Every purely periodic subgraph of $G$ is $(2,2)$-sparse.
		\end{enumerate}
    \end{definition}
       \begin{theorem}\label{ZReflectivelq}
        Let $q\in(1,\infty)\backslash\{2\}$. Let $(\G,\p)$ be a $\Z^2\rtimes\mathcal{C}_s$-regular framework in $(\R^2,\|\cdot\|_q)$ with underlying $\Z^2\rtimes\mathcal{C}_s$-gain graph $(G,m)$. Then $(\G,\p)$ is minimally $\Z^2\rtimes\mathcal{C}_s$-symmetrically infinitesimally rigid if and only if $(G,m)$ is $(\Z^2\rtimes\mathcal{C}_s)_q$-tight.
    \end{theorem}
    As with the other cases, the proof of necessity here is straightforward. Note that there is a $1$-dimensional space of trivial $\Z^2\rtimes\mathcal{C}_s$-symmetric infinitesimal motions (induced by translations parallel to the axis of reflection), so the required rank is $2|V(G)|-1$. We focus on proving sufficiency, and again we will do this by using an inductive construction.
    The relevant extension types are gained $0$-extensions, $1$-extensions, loop-$1$-extensions, vertex-to-$4$-cycle moves and vertex-to-$K_4$ moves. While the $1$-extension requires extra work to show that it preserves minimal $\Z^2\rtimes\mathcal{C}_s$-symmetric rigidity, the other extensions follow readily from standard adaptations.
    \begin{proposition}\label{LqZCsMostExtensions}
		Let $(\G,\p)$ be a $\Z^2\rtimes\mathcal{C}_s$-symmetric framework in $(\R^2,\|\cdot\|_q)$ that is minimally $\Z^2\rtimes\mathcal{C}_s$-symmetrically infinitesimally rigid. Let $(G,m)$ be the underlying $\Z^2\rtimes\mathcal{C}_s$-gain graph of $(\G,\p)$. Let $(G',m')$ be formed by a gained $0$-extension, loop-$1$-extension, vertex-to-$4$-cycle move or vertex-to-$K_4$ move of $(G,m)$. Let $(\G',\p')$ be a $\Z^2\rtimes\mathcal{C}_s$-regular framework derived from $(G',m')$. Then $(\G',\p')$ is minimally $\Z^2\rtimes\mathcal{C}_s$-symmetrically infinitesimally rigid. 
	\end{proposition}
    \begin{proof}
        The proofs that gained $0$-extensions, loop-$1$-extensions and vertex-to-$4$-cycle moves preserve minimal $\Z^2\rtimes\mathcal{C}_s$-symmetric rigidity follow by adapting the methods seen in \cite[Lemma 3.8]{NonEuclidean}, \cite[Lemma 6.1]{EGRES} and \cite[Lemma 3.12]{NonEuclidean} respectively. The proof for gained vertex-to-$K_4$ moves follows by adapting the proof of Proposition \ref{LqCsK4Extension}. See \cite{EssonThesis} for details.
    \end{proof}
    \begin{proposition}\label{LqZCs1Extension}
		Let $(\G,\p)$ be a $\Z^2\rtimes\mathcal{C}_s$-symmetric framework in $(\R^2,\|\cdot\|_q)$ that is minimally $\Z^2\rtimes\mathcal{C}_s$-symmetrically infinitesimally rigid. Let $(G,m)$ be the underlying $\Z^2\rtimes\mathcal{C}_s$-gain graph of $(\G,\p)$. Let $(G',m')$ be formed by a gained $1$-extension of $(G,m)$. Let $(\G',\p')$ be a $\Z^2\rtimes\mathcal{C}_s$-regular framework derived from $(G',m')$. Then $(\G',\p')$ is minimally $\Z^2\rtimes\mathcal{C}_s$-symmetrically infinitesimally rigid.
	\end{proposition}
    \begin{proof}
        In the cases where the new vertex in the gain graph has $2$ or $3$ neighbours, the usual method can be applied. See \cite[Lemma 6.1]{EGRES}, for example. This involves using regularity to assume that the derived neighbours of the new vertex are not collinear.

        This leaves the case where a $1$-extension is performed on a loop to give a triple of parallel edges, as illustrated in Figure \ref{Fig1ExtensionLoop1}. A different method is required to complete the proof in this case, which we now describe. The method is similar to that seen in \cite[Proposition 4.3]{EssonCrystallographic}.

        Suppose that $(G',m')$ is formed from $(G,m)$ by a gained $1$-extension that removes a loop $e=(v_1,v_1;m(e))\in E(G)$ and adds a new vertex $v_0$ with incident edges $e_1=(v_0,v_1;m'(e_1))$, $e_2=(v_0,v_1;m'(e_2))$ and $e_3=(v_0,v_1;m'(e_3))$ such that $(m'(e_1))^{-1}m'(e_2)=m(e)$. A switching operation involving a reflection allows us to assume that $m'(e_3)$ has a trivial $\mathcal{C}_s$-component. A switching operation with a horizontal translation allows us to assume that $m'(e_1)$ has a trivial horizontal translation component and a switching operation with a vertical translation allows us to assume that $m'(e_2)$ has a trivial vertical translation component. Since $(\G,\p)$ is minimally $\Z^2\rtimes\mathcal{C}_s$-symmetrically infinitesimally rigid, every row of the corresponding orbit rigidity matrix must be non-zero. Since a loop with translation gain corresponds to a zero row in the matrix, it can be seen that $m(e)$ has a non-trivial $\mathcal{C}_s$-component and therefore exactly one of $m'(e_1)$ or $m'(e_2)$ has a non-trivial $\mathcal{C}_s$-component. Without loss of generality, suppose that this is $m'(e_1)$. Hence, let $m'(e_1)=(0,d_1,s)$, $m'(e_2)=(c_2,0,0)$ and $m'(e_3)=(c_3,d_3,0)$.
        
        Choose a $\Z^2\rtimes\mathcal{C}_s$-regular configuration $p:V(G)\to\R^2$ of $(G,m)$ (deriving $\p$) and extend this to a configuration $p':V(G')\to\R^2$ of $(G',m')$ (deriving $\p'$) with $p'|_{V(G)}=p$.
        
        We now split the proof into two cases: the case where $d_3\neq0$ and the case where $d_3=0$. First, suppose that $d_3\neq0$. For this case, consider a configuration that aligns $p'(v_0)$ horizontally with $p'(v_1)$. Let $p'(v_1)=(a_1,b)$ and $p'(v_0)=(a_0,b)$ for some $a_0,a_1,b\in\R$, where $b\neq 0$. Suppose that some $u:V(G')\to\R^2$ derives a $\Z^2\rtimes\mathcal{C}_s$-symmetric infinitesimal motion $\tilde{u}$ of $(\G',\p')$. By adding a trivial infinitesimal motion induced by a horizontal translation, it can be assumed that $u(v_1)$ is a vertical vector. Hence, let $u(v_1)=(0,y_1)$ and $u(v_0)=(x_0,y_0)$ for some $x_0,y_0,y_1\in\R$. We aim to show that $u(v_0)=u(v_1)=(0,0)$.

        The constraint on $u$ that is imposed by $e_2$ is
		\begin{align*}
			\begin{pmatrix}
				\mathrm{sgn}(a_0-a_1-c_2)|a_0-a_1-c_2|^{q-1}\\
				0
			\end{pmatrix}
			\cdot
			\begin{pmatrix}
				x_0\\
				y_0-y_1
			\end{pmatrix}
			=0.
		\end{align*}
		Equivalently, $x_0(a_0-a_1-c_2)=0$. An appropriate choice of configuration will ensure that $a_0-a_1-c_2\neq0$, so this constraint implies that $x_0=0$. Given this, the constraint on $u$ that is imposed by $e_1$ is
		\begin{align*}
			\begin{pmatrix}
				\mathrm{sgn}(a_0-a_1)|a_0-a_1|^{q-1}\\
				\mathrm{sgn}(2b-d_1)|2b-d_1|^{q-1}
			\end{pmatrix}
			\cdot
			\begin{pmatrix}
				0\\
				y_0+y_1
			\end{pmatrix}
			=0.
		\end{align*}
		Equivalently, $(y_0+y_1)(2b-d_1)=0$. An appropriate choice of $b\in\R\backslash\{0\}$ will ensure that $2b-d_1\neq0$, so this constraint implies that $y_0+y_1=0$. The constraint imposed by $e_3$ is
		\begin{align*}
			\begin{pmatrix}
				\mathrm{sgn}(a_0-a_1-c_3)|a_0-a_1-c_3|^{q-1}\\
				\mathrm{sgn}(-d_3)|-d_3|^{q-1}
			\end{pmatrix}
			\cdot
			\begin{pmatrix}
				0\\
				y_0-y_1
			\end{pmatrix}
			=0.
		\end{align*}
		Equivalently, $d_3(y_0-y_1)=0$. Since $d_3\neq0$, it can be seen that $y_0-y_1=0$. Since it was previously shown that $y_0+y_1=0$, it must be that $y_0=y_1=0$. This shows that $u(v_0)=u(v_1)=(0,0)$. Hence, the infinitesimal motion $\tilde{u}$ is trivial in this case, as otherwise $(\G,\p)$ is $\Z^2\rtimes\mathcal{C}_s$-symmetrically infinitesimally flexible.

        We now consider the case where $d_3=0$, for which a different configuration will be needed. This time, consider aligning the joints vertically with $p'(v_1)=(a,b_1)$ and $p'(v_0)=(a,b_0)$ for some $a,b_0,b_1\in\R$. Again, let $u:V(G')\to\R^2$ with $u(v_1)=(0,y_1)$ and $u(v_0)=(x_0,y_0)$. As before, the aim is to show that $u$ derives a trivial infinitesimal motion $\tilde{u}$ by showing that $u(v_0)=u(v_1)=(0,0)$. The constraint imposed on $u$ by $e_1$ is
		\begin{align*}
			\begin{pmatrix}
				0\\
				\mathrm{sgn}(b_0+b_1-d_1)|b_0+b_1-d_1|^{q-1}
			\end{pmatrix}
			\cdot
			\begin{pmatrix}
				x_0\\
				y_0+y_1
			\end{pmatrix}
			=0.
		\end{align*}
		Equivalently, $(y_0+y_1)(b_0+b_1-d_1)=0$. An appropriate choice of configuration ensures that $b_0+b_1-d_1\neq0$, so this constraint implies that $y_0+y_1=0$. The constraint imposed by $e_2$ is 
		\begin{align*}
			\begin{pmatrix}
				\mathrm{sgn}(-c_2)|-c_2|^{q-1}\\
				\mathrm{sgn}(b_0-b_1)|b_0-b_1|^{q-1}
			\end{pmatrix}
			\cdot
			\begin{pmatrix}
				x_0\\
				y_0-y_1
			\end{pmatrix}
			=0.
		\end{align*}
		The constraint imposed by $e_3$ is
		\begin{align*}
			\begin{pmatrix}
				\mathrm{sgn}(-c_3)|-c_3|^{q-1}\\
				\mathrm{sgn}(b_0-b_1)|b_0-b_1|^{q-1}
			\end{pmatrix}
			\cdot
			\begin{pmatrix}
				x_0\\
				y_0-y_1
			\end{pmatrix}
			=0.
		\end{align*}
		Combining these constraints, it can be seen that
        \begin{align*}
            x_0\mathrm{sgn}(-c_2)|-c_2|^{q-1}=x_0\mathrm{sgn}(-c_3)|-c_3|^{q-1}.
        \end{align*}
        Since parallel edges must have different gains, $c_2\neq c_3$. It therefore follows that $x_0=0$. Applying this to the constraint imposed by $e_2$ shows that
        \begin{align*}
            (y_0-y_1)\mathrm{sgn}(b_0-b_1)|b_0-b_1|^{q-1}=0.
        \end{align*}
        An appropriate choice of configuration ensures that $b_0-b_1\neq0$ and therefore $y_0-y_1=0$. Combining this with earlier findings shows that $y_0=y_1=0$ and thus $u(v_0)=u(v_1)=(0,0)$. Given this, the constraints imposed by edges in $E(G)-e$ imply that $u(v_1)=u(v_2)=...=u(v_{|V(G)|})=(0,0)$. Hence, $\tilde{u}$ is a trivial infinitesimal motion.
    \end{proof}
    The inductive step of the proof of Theorem \ref{ZReflectivelq} is given by the following result.
    \begin{theorem}\label{ZReflectivelqInduction}
        A $\Z^2\rtimes\mathcal{C}_s$-gain graph is $(\Z^2\rtimes\mathcal{C}_s)_q$-tight if and only if it can be constructed from a $(\Z^2\rtimes\mathcal{C}_s)_q$-tight gain graph on $K_1^1$ by a sequence of gained $0$-extensions, $1$-extensions, loop-$1$-extensions, vertex-to-$4$-cycle moves and vertex-to-$K_4$ moves.
    \end{theorem}

    We now present a proof of this theorem by a sequence of several propositions and lemmas. This is similar to inductive constructions of $(2,2,1)$-gain-tight gain graphs. However, it is made different by the fact that $(2,2)$-sparsity is required on all purely periodic subgraphs, not just balanced ones.

    Basic counting arguments show that every $(2,1)$-tight multigraph has a vertex of degree $2$ or $3$. It is therefore enough to show that every vertex of degree $2$ or $3$ in a $(\Z^2\rtimes\mathcal{C}_s)_q$-tight gain graph admits a reduction that preserves the conditions. Since $K_1^1$ is the only $(2,1)$-tight multigraph on a single vertex, any sequence of such reductions will eventually terminate at a gain graph on $K_1^1$.

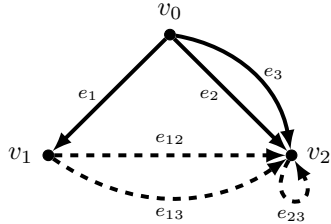
\begin{figure}[H]
        \begin{center}
			\begin{tikzpicture}[scale=1.6]
				\node[vertex,label=above:$v_0$] (v0) at (1,1) {};
				\node[vertex,label=left:$v_1$] (v1) at (0,0) {};
				\node[vertex,label=right:$v_2$] (v2) at (2,0) {};
				
				\draw[dedge] (v0)edge node[left,labelsty]{$e_1$}(v1);
				\draw[dedge] (v0)edge node[left,labelsty]{$e_2$}(v2);
				\draw[dedge] (v0)edge [bend left=35] node[right,labelsty]{$e_3$}(v2);
				
				\draw[dedge,dashed] (v1) to node[above,labelsty]{$e_{12}$}(v2);
				\draw[dedge,dashed] (v1) [bend right=35] to node[below,labelsty]{$e_{13}$}(v2);
				\draw[dedge,dashed] (v2)to [in=-60,out=-120,loop] node[below,labelsty]{$e_{23}$}(v2);
			\end{tikzpicture}
		\end{center}
        \caption{A vertex $v_0$ of degree $3$ with two neighbours.}
        \label{FigDegree3Neighbour2}
    \end{figure}

    Consider each possible neighbourhood of a vertex of degree $2$ or $3$. For a vertex of degree $2$, it is straightforward to see that a $0$-reduction is admissible. 
    \begin{proposition}\label{LqZCsDegree2}
		Suppose that $(G,m)$ is a $(\Z^2\rtimes\mathcal{C}_s)_q$-tight gain graph that has a vertex $v_0$ of degree $2$. Form the gain graph $G-v_0$ by a gained $0$-reduction on $G$ at $v_0$. Then $G-v_0$ is a $(\Z^2\rtimes\mathcal{C}_s)_q$-tight gain graph.
	\end{proposition}
    Likewise, it is easy to see that a loop-$1$-reduction is admissible on any vertex of degree $3$ with an incident loop.
    \begin{proposition}\label{LqZCsDegree3Loop}
		Suppose that $(G,m)$ is a $(\Z^2\rtimes\mathcal{C}_s)_q$-tight gain graph that has a vertex $v_0$ of degree $3$ with an incident loop. Form the gain graph $G-v_0$ by a gained loop-$1$-reduction on $G$ at $v_0$. Then $G-v_0$ is a $(\Z^2\rtimes\mathcal{C}_s)_q$-tight gain graph.
	\end{proposition}
    Now consider the case of a vertex of degree $3$ with one neighbour.
    \begin{proposition}\label{LqZCsDegree3Neighbour1}
		Suppose that $(G,m)$ is a $(\Z^2\rtimes\mathcal{C}_s)_q$-tight gain graph that has a vertex $v_0$ of degree $3$ with a triple of parallel incident edges. Then $(G,m)$ is reducible to a smaller $(\Z^2\rtimes\mathcal{C}_s)_q$-tight gain graph by a gained $1$-reduction at $v_0$.
	\end{proposition}
    \begin{proof}
        Suppose that $v_0$ has a triple of parallel incident edges: $e_1=(v_0,v_1;m(e_1))$, $e_2=(v_0,v_1;m(e_2))$ and $e_3=(v_0,v_1;m(e_3))$ for some $v_1\in V(G)$. To perform a $1$-reduction, remove $v_0$ and add one of the candidate loops: $e_{12}=(v_1,v_1;(m(e_1))^{-1}m(e_2))$, $e_{23}=(v_1,v_1;(m(e_2))^{-1}m(e_3))$ or $e_{31}=(v_1,v_1;(m(e_3))^{-1}m(e_1))$. The aim is therefore to show that one of $G-v_0+e_{12}$, $G-v_0+e_{23}$ or $G-v_0+e_{31}$ is a $(\Z^2\rtimes\mathcal{C}_s)_q$-tight gain graph. Note that there cannot be a loop at $v_1$ in $G-v_0$, as otherwise $\{v_0,v_1\}$ would induce an overcounted subgraph of $G$.

        Since the subgraph induced by $\{v_0,v_1\}$ is $(2,1)$-tight, it is not purely periodic. Hence, at least one of the candidate loops must have a non-trivial $\mathcal{C}_s$-gain component. Without loss of generality, suppose that this is $e_{12}$. With this, consider performing a $1$-reduction to $G-v_0+e_{12}$. Note that any subgraph of $G-v_0+e_{12}$ that contains $e_{12}$ is not purely periodic. Any subgraph of $G-v_0+e_{12}$ that does not contain $e_{12}$ is itself a subgraph of $G$. Hence, $G-v_0+e_{12}$ satisfies the purely periodic condition.

        It is easy to see that any subgraph of $G-v_0+e_{12}$ that contains $e_{12}$ is $(2,1)$-sparse, as otherwise performing the $1$-extension that adds $v_0$ would break the $(2,1)$-sparsity of $G$. Any subgraph of $G-v_0+e_{12}$ that does not contain $e_{12}$ is itself a subgraph of $G$ and is thus $(2,1)$-sparse. Hence, $G-v_0+e_{12}$ satisfies the $(2,1)$-tight condition and is therefore $(\Z^2\rtimes\mathcal{C}_s)_q$-tight.
    \end{proof}
    The next case to consider is a vertex of degree $3$ with two distinct neighbours.
    \begin{proposition}\label{LqZCsDegree3Neighbour2}
		Suppose that $(G,m)$ is a $(\Z^2\rtimes\mathcal{C}_s)_q$-tight gain graph that has a vertex $v_0$ of degree $3$ with two distinct neighbours. Then $(G,m)$ is reducible to a smaller $(\Z^2\rtimes\mathcal{C}_s)_q$-tight gain graph by a gained $1$-reduction at $v_0$.
	\end{proposition}
  There are a number of different cases to consider for the proof of Proposition \ref{LqZCsDegree3Neighbour2}, so we present a proof by means of a sequence of lemmas. To begin this, suppose that $v_0$ has incident edges $e_1=(v_0,v_1;m(e_1))$, $e_2=(v_0,v_2;m(e_2))$ and $e_3=(v_0,v_2;m(e_3))$, for some distinct $v_1,v_2\in V(G)$. After deleting $v_0$, the possible options for edges to add for a $1$-reduction are $e_{12}=(v_1,v_2;(m(e_1))^{-1}m(e_2))$, $e_{13}=(v_1,v_2;(m(e_1))^{-1}m(e_3))$ or $e_{23}=(v_2,v_2;(m(e_2))^{-1}m(e_3))$. To prove Proposition \ref{LqZCsDegree3Neighbour2}, the aim is to show that one of $G-v_0+e_{12}$, $G-v_0+e_{13}$ or $G-v_0+e_{23}$ is a $(\Z^2\rtimes\mathcal{C}_s)_q$-tight gain graph. Figure \ref{FigDegree3Neighbour2} illustrates the neighbourhood of $v_0$, with the candidate edges $e_{12}$, $e_{13}$ and $e_{23}$ represented by dashed lines.

    First, consider the case where $e_2$ and $e_3$ have different $\mathcal{C}_s$-gain components and $G-v_0+e_{23}$ is $(2,1)$-tight.
    \begin{lemma}\label{Degree3Neighbour2Loop}
		Suppose that $(G,m)$ is a $(\Z^2\rtimes\mathcal{C}_s)_q$-tight gain graph that has a vertex $v_0$ of degree $3$, with incident edges $e_1=(v_0,v_1;m(e_1))$, $e_2=(v_0,v_2;m(e_2))$ and $e_3=(v_0,v_2;m(e_3))$ such that $m(e_2)$ and $m(e_3)$ have different $\mathcal{C}_s$-components. If $G-v_0+e_{23}$ is $(2,1)$-tight, then it is a $(\Z^2\rtimes\mathcal{C}_s)_q$-tight gain graph.
	\end{lemma}
	\begin{proof}
		Since $G-v_0+e_{23}$ is $(2,1)$-tight, it is not possible for $e_{23}$ to be parallel to any existing loop. It therefore remains only to check that it satisfies the purely periodic condition. Since $m(e_2)$ and $m(e_3)$ have different $\mathcal{C}_s$-components, $(m(e_2))^{-1}m(e_3)$ has a non-trivial $\mathcal{C}_s$-component. Since $e_{23}$ is a loop,  any subgraph of $G-v_0+e_{23}$ that contains $e_{23}$ will not be purely periodic. Any subgraph of $G-v_0+e_{23}$ that does not contain $e_{23}$ is itself a subgraph of $G$ and therefore satisfies the required conditions on subgraphs. Hence, $G-v_0+e_{23}$ is $(\Z^2\rtimes\mathcal{C}_s)_q$-tight.
	\end{proof}
    For the other cases, the aim is to perform a reduction to either $G-v_0+e_{12}$ or $G-v_0+e_{13}$. Any subgraph of either of these that does not contain the new edge is also a subgraph of $G$ and therefore satisfies all of the conditions required on subgraphs for $(\Z^2\rtimes\mathcal{C}_s)_q$-tightness. Thus it is only necessary to consider subgraphs that contain the new edge. The cases where $G-v_0+e_{12}$ or $G-v_0+e_{13}$ fails each condition can be characterised in terms of subgraphs of $G$, known as blockers. These are described in the following definition.
    \begin{definition}\label{Blockers3}
        Let $(G,m)$ be a $(\Z^2\rtimes\mathcal{C}_s)_q$-tight gain graph. Let $v_0\in V(G)$ be a vertex of degree $3$, with two of its incident edges being $e_i=(v_0,v_i;m(e_i))$ and $e_j=(v_0,v_j;m(e_j))$, for some distinct $v_i,v_j\in V(G)$. Let $e_{ij}=(v_i,v_j;(m(e_i))^{-1}m(e_j))$.
		
		A \emph{$(2,1)$-tight blocker} for the $1$-reduction to $G-v_0+e_{ij}$ is a $(2,1)$-tight subgraph $G_{ij}\subset G$ such that $v_i,v_j\in V(G_{ij})$ and $v_0\notin V(G_{ij})$.

        A \emph{purely periodic blocker} for the $1$-reduction to $G-v_0+e_{ij}$ is a $(2,2)$-tight purely periodic  subgraph $G_{ij}\subset G$ such that $v_i,v_j\in V(G_{ij})$, $v_0\notin V(G_{ij})$ and every path in $G_{ij}$ from $v_i$ to $v_j$ has a net gain with the same $\mathcal{C}_s$-component as $(m(e_i))^{-1}m(e_j)$.
    \end{definition}
    
    Now consider the case where $e_2$ and $e_3$ have different $\mathcal{C}_s$-gain components and $G-v_0+e_{23}$ is  not $(2,1)$-tight.
    \begin{lemma}\label{Degree3Neighbour2OverLoop}
		Suppose that $(G,m)$ is a $(\Z^2\rtimes\mathcal{C}_s)_q$-tight gain graph that has a vertex $v_0$ of degree $3$, with incident edges $e_1=(v_0,v_1;m(e_1))$, $e_2=(v_0,v_2;m(e_2))$ and $e_3=(v_0,v_2;m(e_3))$ such that $m(e_2)$ and $m(e_3)$ have different $\mathcal{C}_s$-components. Suppose that $G-v_0+e_{23}$ is not $(2,1)$-tight. Then one of $G-v_0+e_{12}$ or $G-v_0+e_{13}$ is a $(\Z^2\rtimes\mathcal{C}_s)_q$-tight gain graph.
	\end{lemma}
    \begin{proof}
		Since $G-v_0+e_{23}$ is not $(2,1)$-tight, there is a $(2,1)$-tight subgraph $G_{23}\subseteq G$ with $v_2\in V(G_{23})$ and $v_0\notin V(G_{23})$. Note that $v_1\notin V(G_{23})$, as otherwise adding $v_0$ with its incident edges would break the $(2,1)$-sparsity of $G$. It is easy to see that both $G-v_0+e_{12}$ and $G-v_0+e_{13}$ are $(2,1)$-tight, as otherwise performing the $1$-extension on the new edge would result in an overcounted subgraph of $G$.

        Note that one of $e_{12}$ or $e_{13}$ is not parallel to an existing edge with the same gain, for otherwise $V(G_{23})\cup\{v_0,v_1\}$ contradicts $(2,1)$-sparsity (see Figure~\ref{FigDegree3Neighbour2OverLoop}). So suppose without loss of generality that  $G-v_0+e_{12}$ is a valid gain graph.

        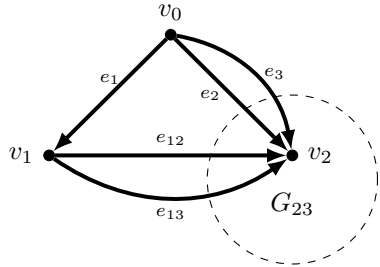
\begin{figure}[H]
            \begin{center}
		\begin{tikzpicture}[scale=1.6]
			\node[vertex,label=above:$v_0$] (v0) at (1,1) {};
			\node[vertex,label=left:$v_1$] (v1) at (0,0) {};
			\node[vertex,label=right:$v_2$] (v2) at (2,0) {};
			
			\draw[dedge] (v0)edge node[above,labelsty]{$e_1$}(v1);
			\draw[dedge] (v0)edge node[left,labelsty]{$e_2$}(v2);
			\draw[dedge] (v0)edge [bend left=35] node[right,labelsty]{$e_3$}(v2);
			
			\draw[edge] (v1) to node[above,labelsty]{}(v2);
            \draw[edge] (v1) to [bend right=35] node[below,labelsty]{}(v2);

            \draw[dashed](2,-0.2)ellipse(0.7 and 0.7);
            \node[blankvertex] (G23) at (2,-0.4) {$G_{23}$};
		\end{tikzpicture}
	\end{center}
    \caption{The subgraph induced by $V(G_{23})\cup\{v_0,v_1\}$ in the case where $e_{12}$ and $e_{13}$ are both parallel to existing edges with equal gains.}
    \label{FigDegree3Neighbour2OverLoop}
        \end{figure}
        
        Suppose that $G-v_0+e_{12}$ fails the purely periodic condition, so there exists a purely periodic blocker  $G_{12}\subset G$. Consider $G_{12}\cup G_{23}$.
        Since $G_{12}$ is $(2,2)$-tight and $G_{23}$ is $(2,1)$-tight, we have
		\begin{align}
			|E(G_{12})\cup E(G_{23})|+|E(G_{12})\cap E(G_{23})|
			=2|V(G_{12})\cup V(G_{23})|+2|V(G_{12})\cap V(G_{23})|-3.\label{EqDegree3Neighbour2}
		\end{align}
		If $E(G_{12})\cap E(G_{23})=\emptyset$, then clearly $|E(G_{12})\cap E(G_{23})|\leq2|V(G_{12})\cap V(G_{23})|-2$. Otherwise, since $G_{12}$ is $(2,2)$-tight, $|E(G_{12})\cap E(G_{23})|\leq2|V(G_{12})\cap V(G_{23})|-2$. In either case, it follows from Equation \eqref{EqDegree3Neighbour2} that $|E(G_{12})\cup E(G_{23})|\geq2|V(G_{12})\cup V(G_{23})|-1$. Adding $v_0$ with its incident edges to $G_{12}\cup G_{23}$ breaks the $(2,1)$-sparsity of $G$, which is a contradiction. Hence, $G-v_0+e_{12}$ satisfies the purely periodic condition.
		
        Hence, one of the $1$-reductions gives a $(\Z^2\rtimes\mathcal{C}_s)_q$-tight gain graph.
	\end{proof}
    We now consider the case where $m(e_2)$ and $m(e_3)$ have the same $\mathcal{C}_s$-component.
    \begin{lemma}\label{Degree3Neighbour2MatchLoop}
		Suppose that $(G,m)$ is a $(\Z^2\rtimes\mathcal{C}_s)_q$-tight gain graph that has a vertex $v_0$ of degree $3$, with incident edges $e_1=(v_0,v_1;m(e_1))$, $e_2=(v_0,v_2;m(e_2))$ and $e_3=(v_0,v_2;m(e_3))$ such that $m(e_2)$ and $m(e_3)$ have the same $\mathcal{C}_s$-component. Then one of $G-v_0+e_{12}$ or $G-v_0+e_{13}$ is a $(\Z^2\rtimes\mathcal{C}_s)_q$-tight gain graph.
	\end{lemma}
    \begin{proof}
        Again, it is clear that both candidate reductions preserve $(2,1)$-tightness. We now show that one of $G-v_0+e_{12}$ or $G-v_0+e_{13}$ is a valid $\Z^2\rtimes\mathcal{C}_s$-gain graph, by checking that there is no pair of parallel edges with equal gains. For contradiction, suppose that both of them fail this, so there is a pair of existing edges in $G$ from $v_1$ to $v_2$ with the same gains as $e_{12}$ and $e_{13}$. This means that $\{v_0,v_1,v_2\}$ induces a $(2,1)$-tight purely periodic subgraph of $G$.
		This contradicts the fact that $(G,m)$ satisfies the purely periodic condition. Hence, one of $e_{12}$ or $e_{13}$ is not parallel to any existing edge of the same gain and therefore the corresponding reduction gives a valid gain graph. Without loss of generality, suppose that $G-v_0+e_{12}$ is a valid $\Z^2\rtimes\mathcal{C}_s$-gain graph. 

        Suppose that $G-v_0+e_{12}$ fails the purely periodic condition, so there is a purely periodic blocker $G_{12}\subseteq G$. Then it can be seen that adding $v_0$ with its incident edges to $G_{12}$ gives a $(2,1)$-tight purely periodic subgraph of $G$, contradicting the fact that $(G,m)$ satisfies the purely periodic condition. This shows that $G-v_0+e_{12}$ is $(\Z^2\rtimes\mathcal{C}_s)_q$-tight.
    \end{proof}
    \begin{proof}[Proof of Proposition \ref{LqZCsDegree3Neighbour2}]
        Combining Lemmas \ref{Degree3Neighbour2Loop}, \ref{Degree3Neighbour2OverLoop} and \ref{Degree3Neighbour2MatchLoop} shows that any vertex of degree $3$ with exactly $2$ neighbours admits a $1$-reduction that preserves the $(\Z^2\rtimes\mathcal{C}_s)_q$-tightness. This completes the proof of Proposition \ref{LqZCsDegree3Neighbour2}.
    \end{proof}

\begin{figure}[H]
        \begin{center}
		\begin{tikzpicture}[scale=1]
			\node[vertex,label=above:$v_0$] (v0) at (1,1) {};
			\node[vertex,label=right:$v_1$] (v1) at (2.5,2) {};
			\node[vertex,label=below:$v_2$] (v2) at (1,-0.5) {};
			\node[vertex,label=left:$v_3$] (v3) at (-0.5,2) {};
			
			\draw[dedge] (v0)edge node[left,labelsty]{$e_1$}(v1);
			\draw[dedge] (v0)edge node[left,labelsty]{$e_2$}(v2);
			\draw[dedge] (v0)edge node[right,labelsty]{$e_3$}(v3);
			
			\draw[dedge,dashed] (v1) to [bend left=30] node[right,labelsty]{$e_{12}$}(v2);
			\draw[dedge,dashed] (v2) to [bend left=30] node[left,labelsty]{$e_{23}$}(v3);
			\draw[dedge,dashed] (v3) to [bend left=30] node[above,labelsty]{$e_{31}$}(v1);
		\end{tikzpicture}
	\end{center}
    \caption{A vertex $v_0$ of degree $3$ with three neighbours.}
    \label{FigDegree3Neighbour3}
    \end{figure}
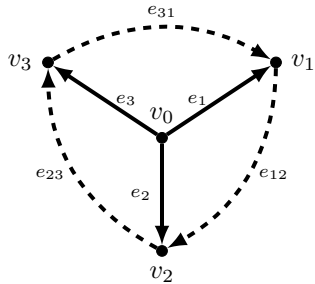

    The last neighbourhood type to consider is a degree $3$ vertex with $3$ distinct neighbours. In this case, the usual reduction to use is a gained $1$-reduction. However, a different reduction will be needed for the specific case when the vertex of degree $3$ is contained in a balanced copy of $K_4$. We will first assume that the vertex of degree $3$ is not in a balanced copy of $K_4$.
    \begin{proposition}\label{Degree3Neighbour3}
		Suppose that $(G,m)$ is a $(\Z^2\rtimes\mathcal{C}_s)_q$-tight gain graph that has a vertex $v_0$ of degree $3$ with $3$ distinct neighbours, which is not contained in a balanced copy of $K_4$. Then $(G,m)$ is reducible to a smaller $(\Z^2\rtimes\mathcal{C}_s)_q$-tight gain graph by a gained $1$-reduction at $v_0$.
	\end{proposition}
    Again, we split the proof of this proposition into a number of lemmas that cover different cases. Suppose that $v_0\in V(G)$ is a vertex of degree $3$ with incident edges $e_1=(v_0,v_1;m(e_1))$, $e_2=(v_0,v_2;m(e_2))$ and $e_3=(v_0,v_3;m(e_3))$, for some distinct $v_1,v_2,v_3\in V(G)$. After deleting $v_0$, the possible options for edges to add for a $1$-reduction are $e_{12}=(v_1,v_2;(m(e_1))^{-1}m(e_2))$, $e_{23}=(v_2,v_3;(m(e_2))^{-1}m(e_3))$ or $e_{31}=(v_3,v_1;(m(e_3))^{-1}m(e_1))$. Figure \ref{FigDegree3Neighbour3} illustrates the neighbourhood of $v_0$, with the candidate edges $e_{12}$, $e_{23}$ and $e_{31}$ represented by dashed lines.

    To prove Proposition \ref{Degree3Neighbour3}, the aim is to show that one of $G-v_0+e_{12}$, $G-v_0+e_{23}$ or $G-v_0+e_{31}$ is a $(\Z^2\rtimes\mathcal{C}_s)_q$-tight gain graph. The following result can be proved by exactly the same proof that was used for \cite[Lemma 5.21]{EssonCrystallographic}.
    \begin{lemma}\label{Condition1Once}
		Suppose that $(G,m)$ is a $(\Z^2\rtimes\mathcal{C}_s)_q$-tight gain graph that has a vertex $v_0$ of degree $3$, which is not in a balanced copy of $K_4$, with edges $e_1$, $e_2$ and $e_3$ from $v_0$ to distinct vertices $v_1$, $v_2$ and $v_3$ respectively. Then at least two of $G-v_0+e_{12}$, $G-v_0+e_{23}$ or $G-v_0+e_{31}$ are $(2,1)$-tight multigraphs.
	\end{lemma}
    Note that Lemma \ref{Condition1Once} does not imply that at least two of $G-v_0+e_{12}$, $G-v_0+e_{23}$ or $G-v_0+e_{31}$ have well-defined gain assignments; this requirement will be considered in subsequent lemmas. We now consider the case where one of the reductions gives a gain graph that is not $(2,1)$-tight.
    \begin{lemma}\label{Fail1Pass2}
		Suppose that $(G,m)$ is a $(\Z^2\rtimes\mathcal{C}_s)_q$-tight gain graph that has a vertex $v_0$ of degree $3$, which is not in a balanced copy of $K_4$, with edges $e_1$, $e_2$ and $e_3$ from $v_0$ to distinct vertices $v_1$, $v_2$ and $v_3$ respectively. Suppose that $G-v_0+e_{12}$ is not $(2,1)$-tight. Then one of $G-v_0+e_{23}$ or $G-v_0+e_{31}$ is a $(\Z^2\rtimes\mathcal{C}_s)_q$-tight gain graph.
	\end{lemma}
    \begin{proof}
		Since $G-v_0+e_{12}$ is not $(2,1)$-tight, there is a $(2,1)$-tight blocker $G_{12}\subseteq G$. We now show that one of $e_{23}$ or $e_{31}$ is not parallel to an existing edge with the same gain. Suppose that both such edges exist. Then $V(G_{12})\cup\{v_0,v_3\}$ induces a $(2,0)$-tight subgraph of $G$, as illustrated in Figure \ref{FigFail1Pass2}. This contradicts the $(2,1)$-sparsity of $G$.
        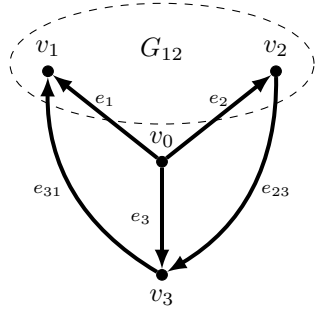
\begin{figure}[H]
            \begin{center}
		\begin{tikzpicture}[scale=1]
			\node[vertex,label=above:$v_0$] (v0) at (1,1) {};
			\node[vertex,label=above:$v_1$] (v1) at (-0.5,2.2) {};
			\node[vertex,label=above:$v_2$] (v2) at (2.5,2.2) {};
			\node[vertex,label=below:$v_3$] (v3) at (1,-0.5) {};
			
			\draw[dedge] (v0)edge node[above,labelsty]{$e_1$}(v1);
			\draw[dedge] (v0)edge node[above,labelsty]{$e_2$}(v2);
			\draw[dedge] (v0)edge node[left,labelsty]{$e_3$}(v3);
			
			\draw[edge] (v2) to [bend left=30] node[right,labelsty]{}(v3);
			\draw[edge] (v3) to [bend left=30] node[left,labelsty]{}(v1);

            \draw[dashed](1,2.3)ellipse(2 and 0.8);
            \node[blankvertex] (G12) at (1,2.5) {$G_{12}$};
		\end{tikzpicture}
	\end{center}
    \caption{The subgraph induced by $V(G_{12})\cup\{v_0,v_3\}$ in the case where $e_{23}$ and $e_{31}$ are both parallel to existing edges with equal gains.}
    \label{FigFail1Pass2}
        \end{figure}
        This shows that one of the reductions gives a valid $\Z^2\rtimes\mathcal{C}_s$-gain graph. Without loss of generality, suppose that this is the reduction to $G-v_0+e_{23}$.
        
        Suppose that $G-v_0+e_{23}$ fails the purely periodic condition, so there is a purely periodic blocker $G_{23}\subset G$. Consider $G_{12}\cup G_{23}$. Since $G_{12}$ is $(2,1)$-tight and $G_{23}$ is $(2,2)$-tight, we have
		\begin{align}
			|E(G_{12})\cup E(G_{23})|+|E(G_{12})\cap E(G_{23})|=2|V(G_{12})\cup V(G_{23})|+2|V(G_{12})\cap V(G_{23})|-3.\label{EqFail1Pass2}
		\end{align}
        If $E(G_{12})\cap E(G_{23})=\emptyset$, then clearly $|E(G_{12})\cap E(G_{23})|\leq2|V(G_{12})\cap V(G_{23})|-2$. Otherwise, since $G_{23}$ is $(2,2)$-tight, $|E(G_{12})\cap E(G_{23})|\leq2|V(G_{12})\cap V(G_{23})|-2$. In either case, it follows from Equation \eqref{EqFail1Pass2} that $|E(G_{12})\cup E(G_{23})|\geq2|V(G_{12})\cup V(G_{23})|-1$. Adding $v_0$ with its incident edges to $G_{12}\cup G_{23}$ breaks $(2,1)$-sparsity of $G$, which is a contradiction. Hence, $G-v_0+e_{23}$ satisfies the purely periodic condition.
	\end{proof}
    We now consider the case where all of the reductions give $(2,1)$-tight multigraphs.
    \begin{lemma}\label{Reflective6}
		Suppose that $(G,m)$ is a $(\Z^2\rtimes\mathcal{C}_s)_q$-tight gain graph that has a vertex $v_0$ of degree $3$, which is not in a balanced copy of $K_4$, with edges $e_1$, $e_2$ and $e_3$ from $v_0$ to distinct vertices $v_1$, $v_2$ and $v_3$ respectively. Suppose that $G-v_0+e_{12}$, $G-v_0+e_{23}$ and $G-v_0+e_{31}$ are all $(2,1)$-tight. Then at least two of them are $(\Z^2\rtimes\mathcal{C}_s)_q$-tight gain graphs.
	\end{lemma}
    \begin{proof}
		Suppose that $G-v_0+e_{12}$ and $G-v_0+e_{23}$ both fail the purely periodic condition. Then they have purely periodic blockers $G_{12}\subset G$ and $G_{23}\subset G$ respectively. By \cite[Lemma 5.17]{EssonCrystallographic},
        \begin{align*}
            |E(G_{12})\cup E(G_{23})|=2|V(G_{12})\cup V(G_{23})|-2\text{ and }|E(G_{12})\cap E(G_{23})|=2|V(G_{12})\cap V(G_{23})|-2.
        \end{align*}
        Since $G_{12}\cap G_{23}$ is $(2,2)$-tight, \cite[Theorem 5]{Pebble} shows that it is connected. Hence, \cite[Lemma 5.20]{EssonCrystallographic} shows that $G_{12}\cup G_{23}$ is purely periodic. Adding $v_0$ with its incident edges to $G_{12}\cup G_{23}$ gives a $(2,1)$-tight purely periodic subgraph of $G$, contradicting the fact that $G$ satisfies the purely periodic condition. This same contradiction can be reached for each pair of candidate edges. Hence, no more than one candidate reduction can give a graph that fails the purely periodic condition.
		
		Now suppose that $e_{12}$ fails the purely periodic condition due to the blocker subgraph $G_{12}$ described above and suppose that $e_{23}$ and $e_{31}$ both have parallel edges with the same gain. Then $V(G_{12})\cup\{v_0,v_3\}$ induces a $(2,1)$-tight balanced subgraph of $G$, contradicting the fact that $G$ satisfies the purely periodic condition.
		
		The only case of failure that has not been ruled out is that where every candidate edge has a parallel edge with the same gain. This means that $v_0$ is in a balanced copy of $K_4$, which contradicts an assumption. Hence, every vertex of degree $3$ that is not in a balanced copy of $K_4$ has an admissible gained $1$-reduction that preserves $(\Z^2\rtimes\mathcal{C}_s)_q$-tightness.
	\end{proof}
    \begin{proof}[Proof of Proposition \ref{Degree3Neighbour3}]
        Combining Lemmas \ref{Condition1Once}, \ref{Fail1Pass2} and \ref{Reflective6} completes the proof of Proposition \ref{Degree3Neighbour3}.
    \end{proof}

    The only case left to consider is a vertex of degree $3$ that is contained within a balanced copy of $K_4$. In this case, there are no admissible $1$-reductions, as any $1$-reduction would result in a pair of parallel edges with equal gains. Instead, the aim is to show that either a $K_4$-to-vertex move or a $4$-cycle-to-vertex move is admissible.

    \begin{proposition}\label{ZReflectivelqBalancedK4}
		Suppose that $(G,m)$ is a $(\Z^2\rtimes\mathcal{C}_s)_q$-tight gain graph that has a vertex $v_0$ of degree $3$, with edges $e_1=(v_0,v_1;m(e_1))$, $e_2=(v_0,v_2;m(e_2))$ and $e_3=(v_0,v_3;m(e_3))$ to distinct $v_1,v_2,v_3\in V(G)$ such that $\{v_0,v_1,v_2,v_3\}$ induces a balanced copy of $K_4$. Then there exists a gained $K_4$-to-vertex move or a gained $4$-cycle-to-vertex move that can be performed on $(G,m)$ to give a smaller $(\Z^2\rtimes\mathcal{C}_s)_q$-tight gain graph.
	\end{proposition}
	\begin{proof}
		By Lemma \ref{lemma:swicth}, it can be assumed that every edge in $K_4(v_0,v_1,v_2,v_3)$ has identity gain. Consider performing a $K_4$-to-vertex move on $K_4(v_0,v_1,v_2,v_3)$. This move always preserves $(\Z^2\rtimes\mathcal{C}_s)_q$-tightness (as it preserves the sparsity and gain space of any subgraph), so the only possibility of failure is by creating a pair of parallel edges with equal gain. This occurs when a pair of vertices in the copy of $K_4$ share a neighbour outside the copy of $K_4$ with edges to the common neighbour that have the same gain. Suppose that this occurs, with $v_1$ and $v_2$ sharing a neighbour $v_4$. By switching operations, it can be assumed that both of the edges to $v_4$ have identity gain. This situation is illustrated in Figure \ref{FigCommonNeighbour}.
		\begin{figure}[H]
		\begin{center}
			\begin{tikzpicture}[scale=1]
				\node[vertex,label=above:$v_0$] (v0) at (1,1) {};
				\node[vertex,label=right:$v_1$] (v1) at (2.5,2) {};
				\node[vertex,label=below:$v_2$] (v2) at (1,-0.5) {};
				\node[vertex,label=left:$v_3$] (v3) at (-0.5,2) {};
				\node[vertex,label=right:$v_4$] (v4) at (3,0.7) {};
				
				\draw[edge] (v0)edge node[left,labelsty]{}(v1);
				\draw[edge] (v0)edge node[left,labelsty]{}(v2);
				\draw[edge] (v0)edge node[right,labelsty]{}(v3);
				
				\draw[edge] (v1) to node[right,labelsty]{}(v2);
				\draw[edge] (v2) to node[left,labelsty]{}(v3);
				\draw[edge] (v3) to node[above,labelsty]{}(v1);
				
				\draw[edge] (v4) to node[right,labelsty]{}(v1);
				\draw[edge] (v4) to node[below,labelsty]{}(v2);
			\end{tikzpicture}
		\end{center}
        \caption{The case where $v_1$ and $v_2$ both have a common neighbour $v_4$, with the edges having identity gain.}
        \label{FigCommonNeighbour}
		\end{figure}
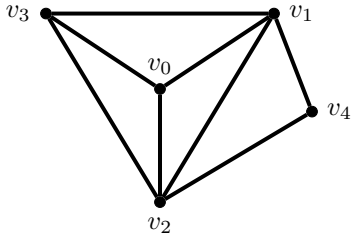
		In this instance, note that $v_0v_1v_4v_2$ forms a $4$-cycle, so consider performing a $4$-cycle-to-vertex move by merging $v_0$ with $v_4$. This particular move is equivalent to removing $v_0$ and adding the edge $(v_3,v_4;(0,0,0))$. The only way that this move could create a pair of parallel edges with equal gain would be if there is an existing edge from $v_3$ to $v_4$ that also has identity gain. However, this would mean that $\{v_0,v_1,v_2,v_3,v_4\}$ induces a $(2,1)$-tight purely periodic subgraph of $G$, contradicting the fact that $G$ satisfies the purely periodic condition. Hence, the $4$-cycle-to-vertex move gives a valid $\Z^2\rtimes\mathcal{C}_s$-gain graph.
		
		Suppose that the $4$-cycle-to-vertex move gives a graph that is not $(2,1)$-tight. This means that there is a $(2,1)$-tight subgraph $G^*\subset G$ with $v_3,v_4\in V(G^*)$ and $v_0\notin V(G^*)$ (as the only new edge added by the move is that from $v_3$ to $v_4$). From here, different contradictions can be reached, depending on whether each of $v_1$ or $v_2$ is in $V(G^*)$.
		\begin{enumerate}
			\item If $v_1,v_2\in V(G^*)$, then adding $v_0$ with its incident edges to $G^*$ breaks $(2,1)$-sparsity of $G$.
			\item If $v_1\in V(G^*)$ and $v_2\notin V(G^*)$ (or vice versa), then adding $v_2$ (or $v_1$) with its incident edges to $G^*$ gives a $(2,0)$-tight subgraph of $G$.
			\item If $v_1,v_2\notin V(G^*)$, then adding $v_1$ and $v_2$ with their incident edges to $G^*$ gives a $(2,0)$-tight subgraph of $G$.
		\end{enumerate}
       		Each case gives a contradiction, so the $4$-cycle-to-vertex move gives a $(2,1)$-tight gain graph.
		
		Suppose that the $4$-cycle-to-vertex move gives a graph that fails the purely periodic condition. This means that there is a $(2,2)$-tight purely periodic subgraph $G^*\subset G$ with $v_3,v_4\in V(G^*)$ and $v_0\notin V(G^*)$ such that every path from $v_3$ to $v_4$ has trivial $\mathcal{C}_s$-gain component. As before, different contradictions can be reached, depending on whether each of $v_1$ or $v_2$ is in $V(G^*)$.
		\begin{enumerate}
			\item If $v_1,v_2\in V(G^*)$, then adding $v_0$ gives a $(2,1)$-tight purely periodic subgraph of $G$. To see that this is purely periodic, note that any $(2,1)$-tight subgraph of $G$ must be vertex-induced and therefore $G^*$ must contain all of the edges induced by $\{v_1,v_2,v_3,v_4\}$. This shows that all paths in $G^*$ between neighbours of $v_0$ must have trivial $\mathcal{C}_s$-gain component and therefore adding $v_0$ preserves the gain space.
			\item If $v_1\in V(G^*)$ and $v_2\notin V(G^*)$ (or vice versa), then adding $v_0$ and $v_2$ (or $v_0$ and $v_1$) with their incident edges to $G^*$ gives a $(2,0)$-tight subgraph of $G$.
			\item If $v_1,v_2\notin V(G^*)$, then adding $v_0$, $v_1$ and $v_2$ with their incident edges to $G^*$ gives a $(2,0)$-tight subgraph of $G$.
		\end{enumerate}
		Each case leads to a contradiction, so the purely periodic condition is satisfied.

	In the case where a $K_4$-to-vertex move fails, it has been shown that there is a $4$-cycle-to-vertex move that gives a $(\Z^2\rtimes\mathcal{C}_s)_q$-tight gain graph. Hence, there is always an admissible reduction.
	\end{proof}
    The proof of Theorem \ref{ZReflectivelqInduction} can now be completed.
    \begin{proof}[Proof of Theorem \ref{ZReflectivelqInduction}]
        It is clear that gained $0$-extensions, $1$-extensions, loop-$1$-extensions, vertex-to-$4$-cycle moves and vertex-to-$K_4$ moves always preserve $(\Z^2\rtimes\mathcal{C}_s)_q$-tightness.
		
		Let $(G,m)$ be a $(\Z^2\rtimes\mathcal{C}_s)_q$-tight gain graph. Basic counting arguments on $(2,1)$-tight graphs show that $G$ has a vertex of degree $2$ or $3$. If $G$ has a vertex of degree $2$, then Proposition \ref{LqZCsDegree2} shows that $(G,m)$ admits a $0$-reduction to a $(\Z^2\rtimes\mathcal{C}_s)_q$-tight gain graph. If $G$ has a vertex of degree $3$ with an incident loop, then Proposition \ref{LqZCsDegree3Loop} shows that $(G,m)$ admits a loop-$1$-reduction to a $(\Z^2\rtimes\mathcal{C}_s)_q$-tight gain graph. If $G$ has a vertex of degree $3$ that has no incident loop and is not in a balanced copy of $K_4$, then, depending on the number of distinct neighbours, Proposition \ref{LqZCsDegree3Neighbour1}, \ref{LqZCsDegree3Neighbour2} or \ref{Degree3Neighbour3} shows that $(G,m)$ admits a $1$-reduction to a $(\Z^2\rtimes\mathcal{C}_s)_q$-tight gain graph. If $G$ has a vertex of degree $3$ that is contained in a balanced copy of $K_4$, then Proposition \ref{ZReflectivelqBalancedK4} shows that $(G,m)$ admits a $K_4$-to-vertex move or $4$-cycle-to-vertex move to a $(\Z^2\rtimes\mathcal{C}_s)_q$-tight gain graph. Hence, every $(\Z^2\rtimes\mathcal{C}_s)_q$-tight gain graph admits a reduction to a smaller $(\Z^2\rtimes\mathcal{C}_s)_q$-tight gain graph.
		
		Since $K_1^1$ is the only $(2,1)$-tight multigraph on a single vertex, repeatedly applying these reductions will eventually lead to a $(\Z^2\rtimes\mathcal{C}_s)_q$-tight gain graph on $K_1^1$. Reverse this sequence of reductions to get the required sequence of extensions to complete the proof.
    \end{proof}
    This allows the proof of sufficiency for Theorem \ref{ZReflectivelq} to be completed.
    \begin{proof}[Proof of sufficiency for Theorem \ref{ZReflectivelq}]
        In any $(\Z^2\rtimes\mathcal{C}_s)_q$-gain graph on $K_1^1$, the loop must have a gain that is not purely periodic. Hence, any $\Z^2\rtimes\mathcal{C}_s$-regular framework derived from a $(\Z^2\rtimes\mathcal{C}_s)_q$-tight gain graph on $K_1^1$ is minimally $\Z^2\rtimes\mathcal{C}_s$-symmetrically infinitesimally rigid. Theorem \ref{ZReflectivelqInduction} states that every $(\Z^2\rtimes\mathcal{C}_s)_q$-tight gain graph can be constructed from a $(\Z^2\rtimes\mathcal{C}_s)_q$-tight gain graph on $K_1^1$ by a sequence of gained $0$-extensions, $1$-extensions, loop-$1$-extensions, vertex-to-$K_4$ moves and vertex-to-$4$-cycle moves. Propositions \ref{LqZCsMostExtensions} and \ref{LqZCs1Extension} show that all of these moves preserve minimal $\Z^2\rtimes\mathcal{C}_s$-symmetric infinitesimal rigidity of $\Z^2\rtimes\mathcal{C}_s$-regular derived frameworks. Hence, every $\Z^2\rtimes\mathcal{C}_s$-regular framework in $(\R^2,\|\cdot\|_q)$ with an underlying gain graph that is $(\Z^2\rtimes\mathcal{C}_s)_q$-tight is minimally $\Z^2\rtimes\mathcal{C}_s$-symmetrically infinitesimally rigid.
    \end{proof}

    \section{Rigidity in the $\ell_1$ and $\ell_\infty$-planes}

    \subsection{Preliminaries on Rigidity and Symmetric Rigidity in the $\ell_1$ and $\ell_\infty$-planes}\label{SectionPolytopicBasics}

    For the rest of this paper, we are interested in frameworks in the $\ell_1$-plane and the $\ell_\infty$-plane. Recall that the \emph{$\ell_1$-norm} $\|\cdot\|_1$ on $\R^2$ is defined, for all $a=(a_1,a_2)\in\R^2$, by
    \begin{align*}
        \|a\|_1=|a_1|+|a_2|.
    \end{align*}
    The \emph{$\ell_\infty$-norm} $\|\cdot\|_\infty$ on $\R^2$ is defined, for all $a=(a_1,a_2)\in\R^2$, by
    \begin{align*}
        \|a\|_\infty=\max\{|a_1|,|a_2|\}.
    \end{align*}
    
    Since they are connected by an isometric isomorphism, we can study just the $\ell_\infty$-plane and the same arguments will also prove characterisations for rigidity in the $\ell_1$-plane and in other polytopic planes with two basis vectors.

    Working with rigidity of frameworks in the $\ell_\infty$-plane requires some methods that differ from those used for other $\ell_q$-planes. These are described in \cite{NonEuclidean} and repeated here for clarity. For each $(x,y)\in\R^2$, define
    \begin{align*}
        \kappa(x,y)=
        \begin{cases}
            (1,0)\text{ if }x>y,\\
            (0,1)\text{ if }x<y,\\
            (0,0)\text{ if }x=y.
        \end{cases}
    \end{align*}

    Let $(G,p)$ be a framework in $(\R^2,\|\cdot\|_\infty)$. As seen in \cite[Proposition 4.2]{NonEuclidean}, an \emph{infinitesimal motion} of $(G,p)$ is a map $u:V(G)\to\R^2$ such that, for all edges $\{v_i,v_j\}\in E(G)$,
    \begin{align*}
        \kappa(p(v_i)-p(v_j))\cdot(u(v_i)-u(v_j))=0.
    \end{align*}

    We can also define the rigidity matrix for a framework in the $\ell_\infty$-plane, which is analogous to the standard rigidity matrix for Euclidean spaces. The \emph{rigidity matrix} $R_{\infty}(G,p)$ of $(G,p)$ is the $|E(G)|\times 2|V(G)|$ matrix with one row for each edge in $E(G)$ and a $2$-tuple of columns for each vertex in $V(G)$, where the row for an edge $e=\{v_i,v_j\}$ is
     \setcounter{MaxMatrixCols}{20}
		\begin{align*}
        \begin{bNiceMatrix}[first-row,first-col]
         &&&&v_i&&&&v_j&&&&\\
		  e&0&...&0&\kappa(p(v_i)-p(v_j))&0&...&0&\kappa(p(v_j)-p(v_i))&0&...&0
        \end{bNiceMatrix}
        .
    \end{align*}
    From this, it is clear that $u\in\R^{2|V(G)|}$ is an infinitesimal motion for $(G,p)$ if and only if $R_{\infty}(G,p)u=0$, so the kernel of the rigidity matrix represents the space of infinitesimal motions. Note that the space of trivial infinitesimal motions is $2$-dimensional, consisting only of infinitesimal motions induced by translations. 

    A framework $(G,p)$ in $(\R^2,\|\cdot\|_\infty)$ is said to be \emph{well-positioned} if, for every edge $\{v_i,v_j\}\in E(G)$, we have $\kappa(p(v_i)-p(v_j))\neq(0,0)$. In this case, an edge $\{v_i,v_j\}\in E(G)$ is said to have \emph{framework colour} $1$ if $\kappa(p(v_i)-p(v_j))=(1,0)$ and $2$ if $\kappa(p(v_i)-p(v_j))=(0,1)$. According to their framework colours, the edges of a well-positioned framework can be partitioned into a pair of \emph{monochrome subgraphs} $G_1$ and $G_2$, consisting of all of the edges of framework colour $1$ and $2$ respectively. This leads to the following useful result.
    \begin{proposition}\label{MonochromeEquivalence}\cite[Propositions 4.3 and 4.4]{NonEuclidean}
        Let $(G,p)$ be a well-positioned framework in $(\R^2,\|\cdot\|_\infty)$. The framework $(G,p)$ is infinitesimally rigid if and only if the monochrome subgraphs $G_1$ and $G_2$ both contain spanning trees of $G$.
    \end{proposition}
    With the assumption that a framework is well-positioned, Proposition \ref{MonochromeEquivalence} provides an alternative way to check infinitesimal rigidity. This is used in \cite{NonEuclidean} to characterise conditions for general rigidity. A variation of this is also used in \cite{GridLike} to characterise conditions for reflectionally-symmetric rigidity and in \cite{KitsonSymmetricNormed} for rotationally-symmetric rigidity. We take a similar approach here to characterise conditions for periodic and $\Z^2\rtimes\mathcal{C}_s$-symmetric rigidity.

    Even when a regular framework is rigid in the $\ell_\infty$ plane, it is important to note that the set of rigid configurations of the underlying graph may not be dense in the set of all configurations, so some care must be taken. However, the set of rigid configurations is always open \cite{DewarDCG21}, so it is possible to choose configurations for convenience to some extent. For example, it can be assumed that no three joints are collinear.

    We now define the orbit rigidity matrix for $\ell_\infty$-planes, which is analogous to that seen for Euclidean planes in \cite{Orbit} and \cite{Inductive}. Since we are only working with symmetries induced by translations and reflections in an axis spanned by either $(1,0)$ or $(0,1)$ (the half-turn symmetry group has been considered in \cite[Theorem 4.3]{KitsonSymmetricNormed} and other symmetry groups will be discussed in Section \ref{SubsectionOtherSymmetries}), our definition will only cover these types of symmetries. Let $\Gamma$ be a group of such isometries on $(\R^2,\|\cdot\|_\infty)$ and let $(G,m)$ be a $\Gamma$-gain graph that derives a $\Gamma$-symmetric framework $(\G,\p)$ in $(\R^2,\|\cdot\|_\infty)$. The \emph{orbit rigidity matrix} $\mathbf{O}_\infty(\G,\p,\Gamma)$ of $(\G,\p)$ is the $|E(G)|\times 2|V(G)|$ matrix with one row for each edge in $(G,m)$ and a $2$-tuple of columns for each vertex in $(G,m)$, defined as follows. For an edge $e=(v_i,v_j;m(e))$ the format of the row for $e$ depends on the type of edge and its gain:
		\begin{enumerate}
			\item Suppose that $e$ is not a loop $(v_i\neq v_j)$. Then the row for $e$ is
            \setcounter{MaxMatrixCols}{20}
		\begin{align*}
        \begin{bNiceMatrix}[first-row,first-col]
         &&&&v_i&&&&v_j&&&&\\
		  e&0&...&0&\kappa(p(v_i)-m(e)p(v_j))&0&...&0&\kappa(p(v_j)-(m(e))^{-1}p(v_i))&0&...&0
        \end{bNiceMatrix}
        .
    \end{align*}
			\item Suppose that $e$ is a loop $(v_i=v_j)$ and the gain of $e$ is a translation. Then the row for $e$ is a zero row.
            \item Suppose that $e$ is a loop $(v_i=v_j)$, the gain of $e$ is a reflection or glide reflection and the framework colour corresponds to a vector that is parallel to the axis of reflection. Then the row for $e$ is a zero row.
            \item Suppose that $e$ is a loop $(v_i=v_j)$, the gain of $e$ is a reflection or glide reflection and the framework colour corresponds to a vector that is perpendicular to the axis of reflection. Then the row for $e$ is
            \begin{align*}
        \begin{bNiceMatrix}[first-row,first-col]
         &&&&v_i&&&&\\
         e&0&...&0&\kappa(p(v_i)-m(e)p(v_i))&0&...&0
        \end{bNiceMatrix}
        .
    \end{align*}
       		\end{enumerate}

    As in other spaces, the orbit rigidity matrix encodes the constraints that are placed on symmetric infinitesimal motions. The space of symmetric infinitesimal motions is given by the kernel of the orbit rigidity matrix. The proof of this is an immediate extension of that seen for the Euclidean orbit rigidity matrix in \cite{Orbit}. See \cite{EssonThesis} for details.

    \subsection{Periodic Rigidity in the $\ell_1$ and $\ell_\infty$-planes}\label{SectionPolytopicPeriodic}
    In this section, we study periodic frameworks in the spaces $(\R^2,\|\cdot\|_1)$ and $(\R^2,\|\cdot\|_\infty)$. As previously mentioned, an isometric isomorphism allows us to focus only on the $\ell_\infty$-plane, as the same arguments will apply directly to the $\ell_1$-plane and other polytopic-normed planes with two basis vectors. Again, we are only studying fixed-lattice periodic rigidity in this paper. Without loss of generality, we assume that the periodicity lattice vectors are $(1,0)$ and $(0,1)$, since the choice of basis vectors is irrelevant in the geometric arguments of the proof. Note that some loops may have no configurations for which they are well-positioned. For example, a loop of gain $(1,1)$ cannot be well-positioned in the $\ell_\infty$-plane with this lattice. However, this is not an issue for characterising minimal rigidity, as loops are always redundant anyway. The following theorem characterises conditions for periodic rigidity in this space.
    
    \begin{theorem}\label{PolyPeriodic}
        Let $(\G,\p)$ be a $\Z^2$-regular framework in $(\R^2,\|\cdot\|_{\infty})$ with underlying $\Z^2$-gain graph $(G,m)$. Then $(\G,\p)$ is minimally periodically infinitesimally rigid if and only if $(G,m)$ is $(2,2)$-tight.
    \end{theorem}
    As with previous results in this paper, necessity follows by straightforward arguments, using the fact  that there is a $2$-dimensional space of trivial  infinitesimal motions.
    
    To prove the sufficiency part of Theorem \ref{PolyPeriodic}, we first need an analogue of Proposition \ref{MonochromeEquivalence} that can be applied to periodic frameworks. Note that the idea of a framework colouring can be applied directly to a periodic framework, with framework colours being constant across each edge orbit, as translations do not change the framework colour of a bar. As such, this can also be thought of as a colouring of the gain graph. Given a framework colouring of $(G,m)$, the \emph{monochrome subgraphs} $G_1$ and $G_2$ are the subgraphs of $G$ consisting of all edges corresponding to edge orbits that receive framework colour $1$ and $2$ respectively. 
    \begin{proposition}\label{MonochromeEquivalencePeriodic}
		Let $(G,m)$ be a $\Z^2$-gain graph that derives a well-positioned periodic framework $(\G,\p)$ in $(\R^2,\|\cdot\|_\infty)$. The framework $(\G,\p)$ is periodically infinitesimally rigid if and only if the corresponding monochrome subgraphs $G_1$ and $G_2$ of $G$ both contain spanning trees of $G$.
	\end{proposition}
    \begin{proof}
        This can be proved in exactly the same way that Proposition \ref{MonochromeEquivalence} was proved in \cite[Propositions 4.3 and 4.4]{NonEuclidean}.
    \end{proof}
    With this, Theorem \ref{PolyPeriodic} can be proved using an inductive method involving gained $0$-extensions, $1$-extensions, edge-to-$K_3$ moves and vertex-to-$K_4$ moves. The first step is to show that each of these moves preserves minimal rigidity of $\Z^2$-gain graphs. The methods for doing this are based on those used for simple graphs in \cite[Lemma 4.9]{NonEuclidean}, showing that each of the extensions preserves the monochrome subgraph property of Proposition \ref{MonochromeEquivalencePeriodic}.
    \begin{proposition}\label{PolyPeriodicExtensions}
		Let $(\G,\p)$ be a $\Z^2$-regular framework in $(\R^2,\|\cdot\|_{\infty})$ that is minimally periodically infinitesimally rigid. Let $(G,m)$ be the underlying $\Z^2$-gain graph of $(\G,\p)$. Let $(G',m')$ be formed by a gained $0$-extension, $1$-extension, edge-to-$K_3$ move or vertex-to-$K_4$ move of $(G,m)$. Let $(\G',\p')$ be a $\Z^2$-regular framework derived from $(G',m')$. Then $(\G',\p')$ is minimally periodically infinitesimally rigid.
	\end{proposition}
	\begin{proof}
        For each type of extension, we use the monochrome subgraph property of Proposition \ref{MonochromeEquivalencePeriodic}. Given $(\G,\p)$, let $p:V(G)\to\R^2$ be the corresponding configuration of the gain graph $(G,m)$. Since $(\G,\p)$ is minimally periodically infinitesimally rigid, Proposition \ref{MonochromeEquivalencePeriodic} shows that the monochrome subgraphs $G_1$ and $G_2$ are both spanning trees of $G$. To show that $\G'$ has a minimally periodically rigid configuration, we aim to find a configuration $p':V(G')\to\R^2$ of $(G',m')$, extending $p$, such that the resulting monochrome subgraphs of $(G',m')$ are both spanning trees. The result will then follow by Proposition \ref{MonochromeEquivalencePeriodic}.
    
		Begin by considering a gained $0$-extension that adds the vertex $v_0$ with edges $e_1=(v_0,v_1;m'(e_1))$ and $e_2=(v_0,v_2;m'(e_2))$, for some (not necessarily distinct) vertices $v_1,v_2\in V(G)$. Consider the lines $L_1=\{m'(e_1)p(v_1)+(t,0):t\in\R\}$ and $L_2=\{m'(e_2)p(v_2)+(0,t):t\in\R\}$. Assuming that the intersection of $L_1$ with $L_2$ is neither coincident to $m'(e_1)p(v_1)$ nor $m'(e_2)p(v_2)$, placing $p'(v_0)$ at this intersection will give framework colour $1$ to $e_1$ and $2$ to $e_2$. By choosing an appropriate configuration $p$ and using the fact that parallel edges must have different gains, it can be seen that $m'(e_1)p(v_1)\neq m'(e_2)p(v_2)$, so no more than one of these can be at the intersection of $L_1$ with $L_2$. If the intersection is at $m'(e_1)p(v_1)$, then it is possible to place $p'(v_0)$ on $L_1$, which gives colour $1$ to $e_1$, sufficiently close to the intersection that $e_2$ receives colour $2$. Similarly, placing $p'(v_0)$ on $L_2$ close to the intersection will achieve this when the intersection is at $m'(e_2)p(v_2)$. Such a placement ensures that both monochrome subgraphs are spanning trees.
		
		Now consider a gained $1$-extension that removes an edge $e=(v_1,v_2;m(e))\in E(G)$ and adds a vertex $v_0$ with edges $e_1=(v_0,v_1;m'(e_1))$, $e_2=(v_0,v_2;m'(e_2))$ and $e_3=(v_0,v_3;m'(e_3))$, for some $v_3\in V(G)$ such that $m(e)=(m'(e_1))^{-1}m'(e_2)$. Without loss of generality, assume that $e$ has framework colour $1$. Since $(\G,\p)$ is minimally periodically rigid, $G$ is $(2,2)$-tight and thus $e$ is not a loop. Hence, $v_1\neq v_2$. A suitable choice of configuration therefore ensures that the points $\{m'(e_1)p(v_1),m'(e_2)p(v_2),m'(e_3)p(v_3)\}$ are not collinear. Let $L_1$ be the line through $m'(e_1)p(v_1)$ and $m'(e_2)p(v_2)$. Let $L_2=\{m'(e_3)p(v_3)+(0,t):t\in\R\}$. Unless this makes it coincident to one of its neighbours, placing $p'(v_0)$ at the intersection of $L_1$ with $L_2$ will give colour $1$ to $e_1$ and $e_2$ and will give colour $2$ to $e_3$. Since the points $\{m'(e_1)p(v_1),m'(e_2)p(v_2),m'(e_3)p(v_3)\}$ are not collinear, the intersection of $L_1$ and $L_2$ is not at $m'(e_3)p(v_3)$. If the intersection is at $m'(e_1)p(v_1)$ or $m'(e_2)p(v_2)$, then placing $p'(v_0)$ on $L_1$ sufficiently close to this intersection will also achieve this framework colouring. Hence, both monochrome subgraphs are spanning trees.
		
		Consider a gained edge-to-$K_3$ move on an edge $e=(v_1,v_2;m(e))\in E(G)$ that adds a vertex $v_0$ with new edges $e_1=(v_0,v_1;m'(e_1))$ and $e_2=(v_0,v_2;m'(e_2))$ such that $m(e)=(m'(e_1))^{-1}m'(e_2)$. Without loss of generality, assume that $e$ has colour $1$. For some $\epsilon>0$, consider setting $p'(v_0)=m'(e_1)p(v_1)+(0,\epsilon)$. For sufficiently small $\epsilon$, $e_1$ receives framework colour $2$ and $e_2$ receives framework colour $1$. Also, all other edges incident to $v_1$ that are moved to be incident to $v_0$ will retain their colour for a sufficiently small $\epsilon$. It then follows that both monochrome subgraphs are spanning trees.
		
		Finally, consider a gained vertex-to-$K_4$ move that replaces a vertex $v_1\in V(G)$ with a balanced $K_4(v_1,w_2,w_3,w_4)$. As seen in \cite[Example 4.5]{NonEuclidean}, the balanced copy of $K_4$ has a configuration in which the monochrome subgraphs are spanning trees. Apply this configuration to $K_4(v_1,w_2,w_3,w_4)$ and rescale it to be sufficiently small that any edges incident to $v_1$ that are moved to be incident to a different vertex in $K_4(v_1,w_2,w_3,w_4)$ retain their colour. Then both monochrome subgraphs of $G'$ are spanning trees.
	\end{proof}
    The following result provides the inductive construction for $(2,2)$-tight $\Z^2$-gain graphs. The approach for proving this is based on those used for symmetry groups of order $2$ in \cite[Theorem 22]{GridLike} and \cite[Theorem 15]{Surfaces}. These have been adapted here for $\Z^2$-gain graphs. Note that this is different to the construction seen in Theorem \ref{PeriodicLqInduction}, as edge-to-$K_3$ moves have been used in place of vertex-to-$4$-cycle moves. This is because the former can be more easily proved to preserve the monochrome subgraph property of Proposition \ref{MonochromeEquivalencePeriodic}.
    \begin{theorem}\label{PolyPeriodicInduction}
		A $\Z^2$-gain graph is $(2,2)$-tight if and only if it can be formed from $K_1$ using a sequence of gained $0$-extensions, $1$-extensions, edge-to-$K_3$ moves and vertex-to-$K_4$ moves.
	\end{theorem}
    \begin{proof}
		Let $(G,m)$ be a $(2,2)$-tight $\Z^2$-gain graph. We aim to find a reduction on $(G,m)$ that gives a smaller $(2,2)$-tight $\Z^2$-gain graph. Basic counting arguments show that $G$ contains a vertex of degree $2$ or $3$. If $G$ contains a vertex $v_0$ of degree $2$, then it is straightforward to show that a $0$-reduction to $G-v_0$ gives a $(2,2)$-tight $\Z^2$-gain graph. If $G$ contains a vertex $v_0$ of degree $3$ that is not in a balanced copy of $K_4$, then the approach seen in \cite[Lemma 4 and Theorem 15]{Surfaces} shows that there is a $1$-reduction on $v_0$ that gives a $(2,2)$-tight $\Z^2$-gain graph. This leaves the case where every vertex of degree $3$ is in a balanced copy of $K_4$. A different approach will be needed for this case.
		
		Suppose that $v_0\in V(G)$ is in the balanced subgraph $A_0=K_4(v_0,v_1,v_2,v_3)$. Since $K_4$-to-vertex moves preserve sparsity, a $K_4$-to-vertex move is admissible here unless two of the vertices in $\{v_1,v_2,v_3\}$ share a neighbour with the same gain. This would prevent a vertex-to-$K_4$ move being applied, as it would create a pair of parallel edges with the same gain. Suppose that this occurs, with $v_4$ being adjacent to both $v_1$ and $v_2$. For this case, we aim to show that the gain graph has an admissible $K_3$-to-edge move. Let $A_1$ be the $(2,2)$-tight subgraph induced by $\{v_0,v_1,v_2,v_3,v_4\}$, which is illustrated in Figure \ref{FigA1}. Note that $A_1$ is balanced, so it can be assumed that all of its edges have identity gain.
        
		If $v_1$ and $v_4$ have edges of equal gain to another common neighbour $v_5\notin V(A_1)$, then let $A_2$ be the $(2,2)$-tight subgraph induced by $\{v_0,v_1,v_2,v_3,v_4,v_5\}$, illustrated in Figure \ref{FigA2}. Again, it can be assumed that all edges have identity gain.

        \begin{figure}[H]
			\centering
            \begin{subfigure}{0.4\textwidth}
			    \begin{tikzpicture}[scale=1]
				\node[vertex,label=above:$v_0$] (v0) at (1,1) {};
				\node[vertex,label=above:$v_1$] (v1) at (2.5,2) {};
				\node[vertex,label=below:$v_2$] (v2) at (1,-0.5) {};
				\node[vertex,label=above:$v_3$] (v3) at (-0.5,2) {};
				\node[vertex,label=right:$v_4$] (v4) at (3,0.7) {};
				
				\draw[edge] (v0)edge node[left,labelsty]{}(v1);
				\draw[edge] (v0)edge node[left,labelsty]{}(v2);
				\draw[edge] (v0)edge node[right,labelsty]{}(v3);
				
				\draw[edge] (v1) to node[right,labelsty]{}(v2);
				\draw[edge] (v2) to node[left,labelsty]{}(v3);
				\draw[edge] (v3) to node[above,labelsty]{}(v1);
				
				\draw[edge] (v4) to node[right,labelsty]{}(v1);
				\draw[edge] (v4) to node[below,labelsty]{}(v2);
			\end{tikzpicture}
				\caption{The subgraph $A_1$.}\label{FigA1}
			\end{subfigure}
			\begin{subfigure}{0.4\textwidth}
				\begin{tikzpicture}[scale=1]
				\node[vertex,label=above:$v_0$] (v0) at (1,1) {};
				\node[vertex,label=above:$v_1$] (v1) at (2.5,2) {};
				\node[vertex,label=below:$v_2$] (v2) at (1,-0.5) {};
				\node[vertex,label=above:$v_3$] (v3) at (-0.5,2) {};
				\node[vertex,label=right:$v_4$] (v4) at (3,0.7) {};
				\node[vertex,label=right:$v_5$] (v5) at (4,1.5) {};
				
				\draw[edge] (v0)edge node[left,labelsty]{}(v1);
				\draw[edge] (v0)edge node[left,labelsty]{}(v2);
				\draw[edge] (v0)edge node[right,labelsty]{}(v3);
				
				\draw[edge] (v1) to node[right,labelsty]{}(v2);
				\draw[edge] (v2) to node[left,labelsty]{}(v3);
				\draw[edge] (v3) to node[above,labelsty]{}(v1);
				
				\draw[edge] (v4) to node[right,labelsty]{}(v1);
				\draw[edge] (v4) to node[below,labelsty]{}(v2);
				
				\draw[edge] (v5) to node[right,labelsty]{}(v1);
				\draw[edge] (v5) to node[below,labelsty]{}(v4);
			\end{tikzpicture}
            \caption{The subgraph $A_2$.}\label{FigA2}
			\end{subfigure}
			\caption{The subgraphs $A_1$ and $A_2$.}
			\label{FigA12}
		\end{figure}
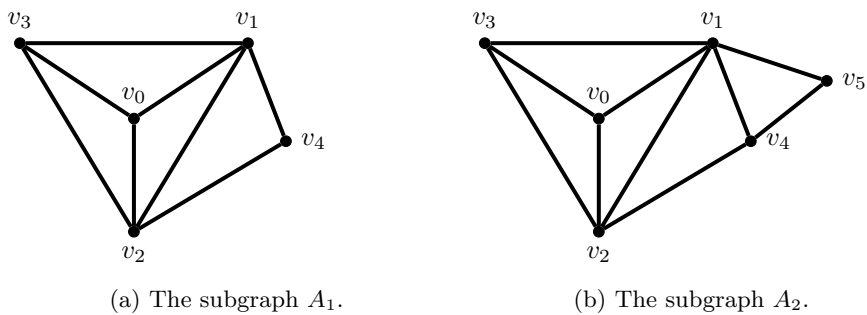
		
		The process of adding more vertices like this will give an increasing sequence $A_0\subset A_1\subset A_2\subset...$, which eventually terminates at some $A_t\subseteq G$. All subgraphs in this sequence are $(2,2)$-tight and balanced. Suppose that $V(A_t\backslash A_{t-1})$ consists of a vertex $v_t$ that is adjacent to $v_i,v_j\in V(A_{t-1})$. Since $A_t$ and $G$ are $(2,2)$-tight, neither $(v_i,v_t;(0,0))$ nor $(v_j,v_t;(0,0))$ is parallel to another edge in $G$. We aim to perform a $K_3$-to-edge move on $\{v_i,v_j,v_t\}$ by contracting $(v_i,v_t;(0,0))$. The lack of a common neighbour of $v_i$ and $v_t$ with the same gain (other than $v_j$) ensures that this will give a valid $\Z^2$-gain graph. We now check that this reduction preserves $(2,2)$-sparsity.
		
		Suppose that the reduced graph is not $(2,2)$-sparse, so there exists a $(2,2)$-tight subgraph $G^*\subseteq G$ such that $\{v_i,v_t\}\in E(G^*)$ and $v_j\notin V(G^*)$. Consider $A_{t-1}\cup G^*$. Since $A_{t-1}$ and $G^*$ are both $(2,2)$-tight, we have
		\begin{align*}
			|E(A_{t-1})\cup E(G^*)|+|E(A_{t-1})\cap E(G^*)|&=2|V(A_{t-1})\cup V(G^*)|+2|V(A_{t-1})\cap V(G^*)|-4.
		\end{align*}
		Since $A_{t-1}$ is $(2,2)$-tight, we have $|E(A_{t-1})\cap E(G^*)|\leq2|V(A_{t-1})\cap V(G^*)|-2$ and thus
			$|E(A_{t-1})\cup E(G^*)|\geq2|V(A_{t-1})\cup V(G^*)|-2$.
		Note that $A_t\cup G^*$ can be obtained from $A_{t-1}\cup G^*$ by merely adding the edge $(v_j,v_t;(0,0))$, so $A_t\cup G^*$ breaks the $(2,2)$-sparsity of $G$. This is a contradiction. Hence, this reduction gives a $(2,2)$-tight $\Z^2$-gain graph.
	\end{proof}
	This completes the proof of sufficiency for Theorem \ref{PolyPeriodic}.
    \begin{proof}[Proof of sufficiency for Theorem \ref{PolyPeriodic}]
        By Theorem \ref{PolyPeriodicInduction}, every $(2,2)$-tight $\Z^2$-gain graph can be formed from $K_1$ using a sequence of gained $0$-extensions, $1$-extensions, edge-to-$K_3$ moves and vertex-to-$K_4$ moves. It is easy to see that any periodic framework derived from $K_1$ is minimally periodically rigid. By Proposition \ref{PolyPeriodicExtensions}, each of the relevant extensions preserves minimal periodic infinitesimal rigidity of $\Z^2$-regular derived frameworks, so every $\Z^2$-regular framework in $(\R^2,\|\cdot\|_\infty)$ derived from a $(2,2)$-tight $\Z^2$-gain graph is minimally periodically rigid.
    \end{proof}

    \subsection{$\Z^2\rtimes\mathcal{C}_s$-symmetric Rigidity in the $\ell_1$ and $\ell_\infty$-Planes}\label{SectionPolytopicWallpaper}
    We now consider $\Z^2\rtimes\mathcal{C}_s$-symmetric rigidity in the $\ell_1$ and $\ell_\infty$-planes, again assuming that the periodicity lattice is fixed. As before, we argue only for the $\ell_\infty$-plane, as the same arguments apply to the $\ell_1$-plane by isometric isomorphism. Similar to \cite{GridLike}, we assume that the axis of reflection in the $\ell_\infty$-plane is spanned by $(1,0)$, which ensures that framework colours are constant over any edge orbit. Recall the definition of $(\Z^2\rtimes\mathcal{C}_s)_q$-tightness from Definition \ref{DefZReflectivelq}.
    \begin{theorem}\label{ZReflectivePoly}
        Let $(\G,\p)$ be a $\Z^2\rtimes\mathcal{C}_s$-regular framework in $(\R^2,\|\cdot\|_{\infty})$ with underlying $\Z^2\rtimes\mathcal{C}_s$-gain graph $(G,m)$. Then $(\G,\p)$ is minimally $\Z^2\rtimes\mathcal{C}_s$-symmetrically infinitesimally rigid if and only if $(G,m)$ is $(\Z^2\rtimes\mathcal{C}_s)_q$-tight.
    \end{theorem}
    Again, necessity can be proved in the standard way, noting that there is a $1$-dimensional space of trivial $\Z^2\rtimes\mathcal{C}_s$-symmetric infinitesimal motions (induced by translations parallel to the axis of reflection). We focus on proving sufficiency. Our method for this combines the approaches of proving other results in this paper, particularly Theorems \ref{ZReflectivelq} and \ref{PolyPeriodic}. First, we need an analogue to the monochrome subgraph property of Proposition \ref{MonochromeEquivalence}. Again, framework colours are constant across each edge orbit, as none of the relevant isometries change the framework colour of a bar. Hence, we can again consider framework colourings on the gain graph. Given a framework colouring of $(G,m)$, the \emph{monochrome subgraphs} $G_1$ and $G_2$ are the subgraphs of $G$ consisting of all edges corresponding to edge orbits that receive framework colour $1$ and $2$ respectively. Recall that a \emph{map graph} is a graph in which every connected component has exactly one cycle (see e.g. \cite{GridLike}).
    \begin{proposition}\label{MonochromeEquivalenceWallpaper}
        Let $(G,m)$ be a $\Z^2\rtimes\mathcal{C}_s$-gain graph that derives a well-positioned $\Z^2\rtimes\mathcal{C}_s$-symmetric framework $(\G,\p)$ in $(\R^2,\|\cdot\|_\infty)$. The framework $(\G,\p)$ is minimally $\Z^2\rtimes\mathcal{C}_s$-symmetrically infinitesimally rigid if and only if the following both hold:
        \begin{enumerate}
            \item The monochrome subgraph $G_1$ is a spanning tree of $G$.
            \item The monochrome subgraph $G_2$ is a spanning map graph of $G$, in which no component is purely periodic.
        \end{enumerate}
    \end{proposition}
    \begin{proof}
        This can be proved in a straightforward way by adapting the method used to prove \cite[Theorem 15]{GridLike}.
    \end{proof}
    As usual, Theorem \ref{ZReflectivePoly} can be proved by an inductive method. The inductive construction uses gained $0$-extensions, $1$-extensions, loop-$1$-extensions, edge-to-$K_3$ moves and vertex-to-$K_4$ moves. Most of the proofs that each of these moves preserves minimal rigidity of a $\Z^2\rtimes\mathcal{C}_s$-regular derived framework follow by adapting standard methods, like those seen in \cite[Theorem 22]{GridLike} and in the proof of Proposition \ref{PolyPeriodicExtensions}. The only cases that require some additional work are those of a $1$-extension on a loop that creates a triple of parallel edges, along with the loop-$1$-extension.
    These were not relevant for $\Z^2$-gain graphs, as $(2,2)$-tight gain graphs never have loops. The $1$-extension to a triple of parallel edges was not relevant for $\mathcal{C}_s$-gain graphs, as triples of parallel edges are not possible with gain graphs for groups of order $2$. On the other hand, loop-$1$-extensions were proved to preserve $\mathcal{C}_s$-symmetric rigidity in \cite[Theorem 22]{GridLike}. However, the proof for $\Z^2\rtimes\mathcal{C}_s$-symmetric rigidity is more complicated, as the translational part of the gain on the loop means that not every position will achieve the required framework colour for the loop.
    \begin{proposition}\label{PolyPeriodicCsExtensions}
		Let $(\G,\p)$ be a $\Z^2\rtimes\mathcal{C}_s$-symmetric framework in $(\R^2,\|\cdot\|_{\infty})$ that is minimally $\Z^2\rtimes\mathcal{C}_s$-symmetrically infinitesimally rigid. Let $(G,m)$ be the underlying $\Z^2\rtimes\mathcal{C}_s$-gain graph of $(\G,\p)$. Let $(G',m')$ be formed by a gained $0$-extension, $1$-extension, loop-$1$-extension, edge-to-$K_3$ move or vertex-to-$K_4$ move of $(G,m)$. Let $(\G',\p')$ be a $\Z^2\rtimes\mathcal{C}_s$-regular framework derived from $(G',m')$. Then $(\G',\p')$ is minimally $\Z^2\rtimes\mathcal{C}_s$-symmetrically infinitesimally rigid.
	\end{proposition}
    \begin{proof}
        Like in the proof of Proposition \ref{PolyPeriodicExtensions}, we can use the monochrome subgraph property of Proposition \ref{MonochromeEquivalenceWallpaper}. Given $(\G,\p)$, let $p:V(G)\to\R^2$ be the corresponding configuration of the gain graph $(G,m)$. Since $(\G,\p)$ is minimally $\Z^2\rtimes\mathcal{C}_s$-symmetrically infinitesimally rigid, Proposition \ref{MonochromeEquivalenceWallpaper} shows that the monochrome subgraph $G_1$ is a spanning tree of $G$ and the monochrome subgraph $G_2$ is a spanning map graph of $G$ in which no component is purely periodic. To show that $\G'$ has a minimally $\Z^2\rtimes\mathcal{C}_s$-symmetrically infinitesimally rigid configuration, we aim to find a configuration $p':V(G')\to\R^2$ of $(G',m')$, extending $p$, such that the resulting monochrome subgraphs of $(G',m')$ satisfy the property of Proposition \ref{MonochromeEquivalenceWallpaper}.
        
        As mentioned earlier, most of the cases for this proof follow from the arguments seen in \cite[Theorem 22]{GridLike}. The only other cases to consider for $\Z^2\rtimes\mathcal{C}_s$-gain graphs are $1$-extensions on loops that create a triple of parallel edges and loop-$1$-extensions.
        
        Suppose that a $1$-extension is performed on the loop $e=(v_1,v_1;m(e))\in E(G)$ by adding the vertex $v_0$ with edges $e_1=(v_0,v_1;m'(e_1))$, $e_2=(v_0,v_1;m'(e_2))$ and $e_3=(v_0,v_1;m'(e_3))$ such that $m(e)=(m'(e_1))^{-1}m'(e_2)$. Since $e$ is a loop, it must have a non-trivial $\mathcal{C}_s$-gain component. By $\Z^2\rtimes\mathcal{C}_s$-regularity, $e$ must have framework colour $2$. It is sufficient to find a configuration of $(G',m')$ that extends $p$ such that $e_1$ and $e_2$ receive framework colour $2$, while $e_3$ receives framework colour $1$. To do this, set $p':V(G')\to\R^2$ such that $p'|_{V(G)}=p$. Let $L_1$ be the line through $m'(e_3)p(v_1)$ in direction $(1,0)$. Let $L_2$ be the line through $m'(e_1)p(v_1)$ and $m'(e_2)p(v_1)$. Then placing $p'(v_0)$ at the intersection of $L_1$ with $L_2$ will give the required framework colours, unless this intersection is coincident to one of $m'(e_1)p(v_1)$, $m'(e_2)p(v_1)$ or $m'(e_3)p(v_1)$. Note that these are all distinct points, so only one of these coincidences can occur. In the case where the intersection is at $m'(e_1)p(v_1)$ or $m'(e_2)p(v_1)$, it is possible to place $p'(v_0)$ on $L_2$, sufficiently close to the intersection to get the required framework colouring. In the case where the intersection is at $m'(e_3)p(v_1)$, it is possible to get the required framework colouring by placing $p'(v_0)$ on $L_1$ sufficiently close to the intersection. Hence, $(\G',\p')$ is minimally $\Z^2\rtimes\mathcal{C}_s$-symmetrically infinitesimally rigid by Proposition \ref{MonochromeEquivalenceWallpaper}.

        It remains only to consider loop-$1$-extensions. Suppose that a loop-$1$-extension adds a vertex $v_0$ with incident edges $l=(v_0,v_0;m'(l))$ and $e=(v_0,v_1;m'(e))$. We aim to choose a position $p'(v_0)$ such that $e$ receives framework colour $1$ and $l$ receives framework colour $2$. For ease of notation, let $p'(v_0)=(a_0,b_0)$, and $p(v_1)=(a_1,b_1)$. Also, let $m'(l)=(c_l,d_l,s)$. By switching operations, it can be assumed that $m'(e)=(0,0,0)$. To ensure that $l$ receives framework colour $2$, it must be that
        \begin{align*}
            \kappa(p'(v_0)-m'(l)p'(v_0))=\kappa(-c_l,2b_0-d_l)=(0,1).
        \end{align*}
        Equivalently,
        \begin{align*}
            |c_l|<|2b_0-d_l|
        \end{align*}
               Clearly, it is possible to choose $b_0\in\R$ such that this holds. To ensure that $e$ receives framework colour $1$, it must be that
        \begin{align*}
            \kappa(p'(v_0)-p(v_1))=\kappa(a_0-a_1,b_0-b_1)=(1,0).
        \end{align*}
        Equivalently,
        \begin{align*}
            |a_0-a_1|>|b_0-b_1|.
        \end{align*}
        Given any values of $b_0,a_1,b_1\in\R$, it is clearly always possible to choose $a_0\in\R$ such that this holds. With the resulting configuration, the framework $(\G',\p')$ satisfies the monochrome subgraph property and is thus minimally $\Z^2\rtimes\mathcal{C}_s$-symmetrically infinitesimally rigid by Proposition \ref{MonochromeEquivalenceWallpaper}.
    \end{proof}

    The following result gives the inductive construction relevant to Theorem \ref{ZReflectivePoly}.
    \begin{theorem}\label{ZReflectivePolyInduction}
        A $\Z^2\rtimes\mathcal{C}_s$-gain graph is $(\Z^2\rtimes\mathcal{C}_s)_q$-tight if and only if it can be constructed from a $(\Z^2\rtimes\mathcal{C}_s)_q$-tight gain graph on $K_1^1$ by a sequence of gained $0$-extensions, $1$-extensions, loop-$1$-extensions, edge-to-$K_3$ moves and vertex-to-$K_4$ moves.
    \end{theorem}
    \begin{proof}
        Basic counting arguments show that every $(2,1)$-tight multigraph has a vertex of degree $2$ or $3$. It is therefore enough to show that every vertex of degree $2$ or $3$ in a $(\Z^2\rtimes\mathcal{C}_s)_q$-tight gain graph admits a reduction that preserves $(\Z^2\rtimes\mathcal{C}_s)_q$-tightness. In the cases of a degree $2$ vertex or a degree $3$ vertex that is not contained in a balanced copy of $K_4$, the methods seen in Section \ref{SectionLqWallpaper} work in exactly the same way. This leaves only the case of a degree $3$ vertex that is contained in a balanced copy of $K_4$. For this case, adapting the method used in the proof of Theorem \ref{PolyPeriodicInduction} shows that there is a $K_4$-to-vertex move or $K_3$-to-edge move that preserves $(\Z^2\rtimes\mathcal{C}_s)_q$-tightness.
    \end{proof}
    This completes the proof of sufficiency for Theorem \ref{ZReflectivePoly}.
     \begin{proof}[Proof of sufficiency for Theorem \ref{ZReflectivePoly}]
        By Theorem \ref{ZReflectivePolyInduction}, any $(\Z^2\rtimes\mathcal{C}_s)_q$-tight gain graph can be formed from a $(\Z^2\rtimes\mathcal{C}_s)_q$-tight gain graph on $K_1^1$ by a sequence of gained $0$-extensions, $1$-extensions, edge-to-$K_3$ moves and vertex-to-$K_4$ moves. Given a $(\Z^2\rtimes\mathcal{C}_s)_q$-tight gain graph on $K_1^1$, it is always possible to choose a configuration where the loop receives framework colour 2. By Proposition \ref{MonochromeEquivalenceWallpaper}, the resulting derived framework is minimally $\Z^2\rtimes\mathcal{C}_s$-symmetrically infinitesimally rigid. By Proposition \ref{PolyPeriodicCsExtensions}, each of the relevant extensions preserves minimal $\Z^2\rtimes\mathcal{C}_s$-symmetric infinitesimal rigidity of $\Z^2\rtimes\mathcal{C}_s$-regular derived frameworks, so every $\Z^2\rtimes\mathcal{C}_s$-regular framework in $(\R^2,\|\cdot\|_\infty)$ derived from a $(\Z^2\rtimes\mathcal{C}_s)_q$-tight gain graph is minimally $\Z^2\rtimes\mathcal{C}_s$-symmetrically infinitesimally rigid.
    \end{proof}

    \section{Further Work}\label{SectionFurther}
    \subsection{Other Symmetries}\label{SubsectionOtherSymmetries}
    It is natural to consider how to characterise forced-symmetric rigidity in non-Euclidean $\ell_q$-planes and polytopic planes with respect to other symmetry groups. The most obvious class of groups would be those involving rotational symmetries. For polytopic planes, conditions for $2$-fold rotationally-symmetric rigidity were characterised by D. Kitson, A. Nixon and B. Schulze in \cite[Theorem 4.3]{KitsonSymmetricNormed}. It still remains open to do this for $4$-fold rotational symmetry or for any rotational symmetry in $\ell_q$-planes for $q\in(1,\infty)\backslash\{2\}$ (where only $2$-fold and $4$-fold rotations can be considered). Characterising forced-symmetric rigidity for these by an inductive method would likely be very challenging, as gain graphs for minimally forced-symmetrically rigid frameworks are $(2,0)$-tight and inductive constructions of $(2,0)$-tight gain graphs are generally quite difficult. However, an inductive construction was successfully obtained for $2$-fold rotationally-symmetric rigidity in polytopic planes in \cite[Theorem 4.3]{KitsonSymmetricNormed}, so a similar construction may work for other cases. In the case of $4$-fold rotational symmetry in polytopic planes, another issue is that framework colours may not be consistent across edge orbits, so there is no clear analogue to Proposition \ref{MonochromeEquivalence}.

    When studying reflectionally-symmetric rigidity in polytopic planes, both this paper and \cite{GridLike} have assumed that the axis of reflection is spanned by one of the vectors used to generate the polytopic norm. It would also be possible to consider reflectional symmetry in a diagonal of the unit circle. This seems like a harder problem, as framework colours will not be consistent across edge orbits and therefore the idea of Proposition \ref{MonochromeEquivalence} cannot be directly applied.

    \subsection{Flexible Lattice Representations}
    When working with periodic frameworks in this paper, we have only been considering infinitesimal motions that are periodic on the fixed lattice. An immediate question is how this could be modified to fully flexible lattice representations or at least to partially flexible lattice representations. For the Euclidean plane, the rigidity of periodic frameworks on a fully flexible lattice was first studied by Borcea and Streinu \cite{borstreiproc}, and a characterisation of periodic rigidity in this setting was obtained by Malestein and Theran in \cite[Theorem A]{MalesteinFlexible}, although this did not use an inductive construction. For a partially-flexible lattice in the Euclidean plane, an inductive construction for minimally rigid gain graphs was found in \cite[Theorem 2]{InductiveFlexible}. Moreover, there are relative characterisations of periodic rigidity (in the case where the gain assignments are not part of the initial data) in all dimensions \cite{BorceaStreinuMinimally,Quotient}.    
    It may be possible to adapt one of these approaches to characterise flexible-lattice periodic rigidity in non-Euclidean planes.

    \subsection{Higher Dimensions}
    Another idea is to study forced-symmetric and forced-periodic rigidity in higher-dimensional $\ell_q$-spaces and polytopic spaces. As in the Euclidean case, this is likely to be very difficult, as there is currently no known inductive construction for the relevant counts. Indeed, there is no known characterisation of non-symmetric rigidity of finite frameworks in any higher-dimensional $\ell_q$-spaces or polytopic spaces. Note that minimally rigid graphs in $3$-dimensional non-Euclidean $\ell_q$-spaces and polytopic spaces must be at least $(3,3)$-tight. An inductive construction for this is not known, but it may be easier than the $(3,6)$-tight construction that would be needed for the Euclidean case.
    
    For body-bar frameworks, a complete characterisation of forced-symmetric rigidity for all symmetry groups in all dimensions of Euclidean spaces was proved by S. Tanigawa in \cite[Theorem 7.2]{TanigawaBodyBar}. A similar proof may be possible for non-Euclidean spaces.

	\bibliographystyle{plain}

\begin{thebibliography}{10}

\bibitem{BernsteinRotations}
D.~I. Bernstein.
\newblock Generic symmetry-forced infinitesimal rigidity: translations and rotations.
\newblock {\em SIAM J. Appl. Algebra Geom.}, 6(2):190--215, 2022.

\bibitem{borstreiproc}
C.~S. Borcea and I.~Streinu.
\newblock Periodic frameworks and flexibility.
\newblock {\em Proc. R. Soc. Lond. Ser. A Math. Phys. Eng. Sci.}, 466(2121):2633--2649, 2010.

\bibitem{BorceaStreinuMinimally}
C.~S. Borcea and I.~Streinu.
\newblock Minimally rigid periodic graphs.
\newblock {\em Bull. Lond. Math. Soc.}, 43(6):1093--1103, 2011.

\bibitem{connelly_guest_2022}
R.~Connelly and S.~D. Guest.
\newblock {\em Frameworks, Tensegrities, and Symmetry}.
\newblock Cambridge University Press, 2022.

\bibitem{Sym}
J.~H. Conway, H.~Burgiel, and C.~Goodman-Strauss.
\newblock {\em The symmetries of things}.
\newblock A K Peters, Ltd., Wellesley, MA, 2008.

\bibitem{Dewar}
S.~Dewar.
\newblock Infinitesimal rigidity in normed planes.
\newblock {\em SIAM J. Discrete Math.}, 34(2):1205--1231, 2020.

\bibitem{DewarDCG21}
S.~Dewar.
\newblock Equivalence of continuous, local and infinitesimal rigidity in normed spaces.
\newblock {\em Discrete Comput. Geom.}, 65(3):655--679, 2021.

\bibitem{Quotient}
S.~Dewar, G.~Grasegger, E.~Kastis, and A.~Nixon.
\newblock Quotient graphs of symmetrically rigid frameworks.
\newblock {\em Doc. Math.}, 29(3):561--595, 2024.

\bibitem{dhn24}
S.~Dewar, J.~Hewetson, and A.~Nixon.
\newblock Coincident-point rigidity in normed planes.
\newblock {\em Ars Math. Contemp.}, 24(1):Paper No. 10, 17, 2024.

\bibitem{DewarCylindrical}
S.~Dewar and D.~Kitson.
\newblock Rigid graphs in cylindrical normed spaces.
\newblock {\em SIAM Journal on Discrete Mathematics}, 39(3):1545--1567, 2025.

\bibitem{EssonThesis}
J.~Esson.
\newblock {\em Crystallographic frameworks in Euclidean and non-Euclidean spaces}.
\newblock PhD thesis, Lancaster University, Expected in 2027.

\bibitem{EssonCrystallographic}
J.~Esson, E.~Kastis, and B.~Schulze.
\newblock Orientation-reversing crystallographic rigidity.
\newblock {\em arXiv:2502.14392}, 2025.

\bibitem{EGRES}
T.~Jord\'an, V.~E. Kaszanitzky, and S.~Tanigawa.
\newblock Gain-sparsity and symmetry-forced rigidity in the plane.
\newblock {\em Discrete Comput. Geom.}, 55(2):314--372, 2016.

\bibitem{KastisSymbol}
E.~Kastis, D.~Kitson, and J.~E. McCarthy.
\newblock Symbol functions for symmetric frameworks.
\newblock {\em J. Math. Anal. Appl.}, 497(2):Paper No. 124895, 24, 2021.

\bibitem{KitsonDenseOpen}
D.~Kitson.
\newblock Finite and infinitesimal rigidity with polyhedral norms.
\newblock {\em Discrete Computational Geometry}, 54:390--411, 2015.

\bibitem{KitsonSymmetricNormed}
D.~Kitson, A.~Nixon, and B.~Schulze.
\newblock Rigidity of symmetric frameworks in normed spaces.
\newblock {\em Linear Algebra Appl.}, 607:231--285, 2020.

\bibitem{NonEuclidean}
D.~Kitson and S.~C. Power.
\newblock Infinitesimal rigidity for non-{E}uclidean bar-joint frameworks.
\newblock {\em Bull. Lond. Math. Soc.}, 46(4):685--697, 2014.

\bibitem{GridLike}
D.~Kitson and B.~Schulze.
\newblock Motions of grid-like reflection frameworks.
\newblock {\em J. Symbolic Comput.}, 88:47--66, 2018.

\bibitem{Laman}
G.~Laman.
\newblock On graphs and rigidity of plane skeletal structures.
\newblock {\em J. Engrg. Math.}, 4:331--340, 1970.

\bibitem{Pebble}
A.~Lee and I.~Streinu.
\newblock Pebble game algorithms and sparse graphs.
\newblock {\em Discrete Math.}, 308(8):1425--1437, 2008.

\bibitem{MalesteinFlexible}
J.~Malestein and L.~Theran.
\newblock Generic combinatorial rigidity of periodic frameworks.
\newblock {\em Adv. Math.}, 233:291--331, 2013.

\bibitem{MalesteinSymmetry}
J.~Malestein and L.~Theran.
\newblock Frameworks with forced symmetry {I}: reflections and rotations.
\newblock {\em Discrete Comput. Geom.}, 54(2):339--367, 2015.

\bibitem{InductiveFlexible}
A.~Nixon and E.~Ross.
\newblock Periodic rigidity on a variable torus using inductive constructions.
\newblock {\em Electron. J. Combin.}, 22(1):Paper 1.1, 30, 2015.

\bibitem{Surfaces}
A.~Nixon and B.~Schulze.
\newblock Symmetry-forced rigidity of frameworks on surfaces.
\newblock {\em Geom. Dedicata}, 182:163--201, 2016.

\bibitem{Pollaczek-Geiringer}
H.~Pollaczek-Geiringer.
\newblock Über die {G}liederung ebener {F}achwerke.
\newblock {\em ZAMM - Journal of Applied Mathematics and Mechanics / Zeitschrift für Angewandte Mathematik und Mechanik}, 7(1):58--72, 1927.

\bibitem{Torus}
E.~Ross.
\newblock The rigidity of periodic frameworks as graphs on a fixed torus.
\newblock {\em Contrib. Discrete Math.}, 9(1):11--45, 2014.

\bibitem{Inductive}
E.~Ross.
\newblock Inductive constructions for frameworks on a two-dimensional fixed torus.
\newblock {\em Discrete Comput. Geom.}, 54(1):78--109, 2015.

\bibitem{HMbook}
D.~E. Sands.
\newblock {\em Introduction to Crystallography}.
\newblock Dover Publications Inc, 2003.

\bibitem{Schatt}
D.~Schattschneider.
\newblock The plane symmetry groups: Their recognition and notation.
\newblock {\em The American Mathematical Monthly}, 85(6):439--450, 1978.

\bibitem{constraint19}
B.~Schulze.
\newblock Combinatorial rigidity of symmetric and periodic frameworks.
\newblock In M.~Sitharam, A.~St.~John, and J.~Sidman, editors, {\em {H}andbook of {G}eometric {C}onstraint {S}ystems {P}rinciples}. CRC Press, 2019.

\bibitem{Orbit}
B.~Schulze and W.~Whiteley.
\newblock The orbit rigidity matrix of a symmetric framework.
\newblock {\em Discrete Comput. Geom.}, 46(3):561--598, 2011.

\bibitem{HandbookDCG}
B.~Schulze and W.~Whiteley.
\newblock Rigidity and scene analysis.
\newblock In {\em Handbook of Discrete and Computational Geometry, Third Edition}, C.D. Toth, J. O'Rourke, J.E. Goodman, editors. Chapman $\&$ Hall CRC, 2018.

\bibitem{HandbookDCGsym}
B.~Schulze and W.~Whiteley.
\newblock Rigidity of symmetric frameworks.
\newblock In {\em Handbook of Discrete and Computational Geometry, Third Edition}, C.D. Toth, J. O'Rourke, J.E. Goodman, editors. Chapman $\&$ Hall CRC, 2018.

\bibitem{Gocg}
A.~Szczepanski.
\newblock {\em Geometry of crystallographic groups}, volume~4 of {\em Algebra and discrete mathematics}.
\newblock World Scientific Publishing, 2012.

\bibitem{TanigawaBodyBar}
S.~Tanigawa.
\newblock Matroids of gain graphs in applied discrete geometry.
\newblock {\em Trans. Amer. Math. Soc.}, 367(12):8597--8641, 2015.

\bibitem{WhiteleyMatroids}
W.~Whiteley.
\newblock Some matroids from discrete applied geometry.
\newblock In {\em Matroid theory ({S}eattle, {WA}, 1995)}, volume 197 of {\em Contemp. Math.}, pages 171--311. Amer. Math. Soc., Providence, RI, 1996.

\end{thebibliography}

\end{document}